\documentclass[a4paper,11pt]{article}
\usepackage{color,xcolor,ucs}
\usepackage[top=0.1in, bottom=0.43in, left = 0.52in, right = 0.52in]{geometry}
\usepackage[linkcolor=black,colorlinks=true,urlcolor=blue]{hyperref}
\usepackage{mathtools}

\usepackage{color,xcolor,ucs}
\usepackage{mathtools}   \usepackage{tikz} 
\usepackage{ amssymb }
\usepackage{extarrows} 
\usepackage{pgf,tikz}
\usepackage{float}
\usetikzlibrary{positioning}
\usetikzlibrary{shapes.geometric}
\usetikzlibrary{shapes.misc}
\usetikzlibrary{arrows}
\usepackage{caption}
\usepackage{mathrsfs}
\usetikzlibrary{arrows,shapes,automata,backgrounds,petri,positioning}
\usetikzlibrary{decorations.pathmorphing}
\usetikzlibrary{decorations.shapes}
\usetikzlibrary{decorations.text}
\usetikzlibrary{decorations.fractals}
\usetikzlibrary{decorations.footprints}
\usetikzlibrary{shadows}
\usetikzlibrary{calc}
\usetikzlibrary{spy}
\usepackage{amsmath}
\usepackage{array}
\usepackage{ amssymb }
\usepackage{braket}
\usepackage{qcircuit}
\usepackage{soul}
\usepackage{braket}

\usepackage{amsmath}
\usepackage{ amssymb }
\usepackage{braket}
\usepackage{qcircuit}
\usepackage{soul}
\usepackage{braket}

\title{{\LARGE Kesten's incipient infinite cluster for the three-dimensional, metric-graph Gaussian free field, from critical level-set percolation, and for the Villain model}}
\author{Pete Rigas \footnote{Cornell AEP and Mathematics Departments, pbr43@cornell.edu, ORCID: 0009-0003-1053-9720}}
\date{}

\begin{document}

\maketitle

\begin{abstract}
    We address one open problem in a recent work due to Ding and Wirth, the first version of which was available in $2019$, relating to level-set percolation on metric-graphs for the Gaussian free field in three dimensions, in which it was shown that a percolation estimate that the authors employ for studying connectivity properties of different heights of the metric graph Gaussian free field is bounded above poly-logarithmically. In three dimensions, in order to construct Kesten's incipient infinite cluster which was first seminally introduced for Bernoulli percolation in two dimensions, in $1986$, through the equality of two probabilistic quantities, we make use of a streamlined version of the $1986$ argument due to Basu and Sapozhnikov, which was first made available in $2016$, that introduces properties of crossing probabilities for demonstrating that the IIC exists for Bernoulli percolation on an infinite connected, bounded degree graph, namely that crossing probabilities are quasi-multiplicative. To make use of such arguments for demonstrating the existence of the metric-graph GFF IIC in three dimensions, we also address another open problem raised in a recent work, from October $2022$, due to Dubedat and Falconet, which expresses an open problem pertaining to the construction of an IIC-type limit for the Villain model. 
\end{abstract}

\section{Introduction}

\subsection{Overview}

The Gaussian Free Field (GFF) is a random surface model that has attracted significant attention within the mathematical community, with efforts dedicated towards studying level-set connectivity properties through construction of percolation interfaces {\color{blue}[1,2,9,12,27,29,30]}, quantification of extrema {\color{blue}[6,8]}, connections to eigenvalues of the corners of random matrices {\color{blue}[7]}, the existence of a phase transition for percolation by constructing GFF level-set dependent copies of edges on graphs whose randomness is integrated out with an Edwards-Sokal coupled measure {\color{blue}[14]}, asymptotics through a renormalization scheme {\color{blue}[15]}, and characterizations of the field through the domain Markov property permitting for a harmonic decomposition over boundaries of finite volumes {\color{blue}[3,11,28]}. In addition to such results, other well known properties that are relevant for open questions raised in {\color{blue}[11]} pertain to the behavior and organization of critical percolation clusters for $d=3$ on metric graphs, which are analyzed with exploration martingales, include previous works providing results on the chemical distance {\color{blue}[10,20]}, sign cluster arrangements {\color{blue}[17,18]}, percolation critical exponents {\color{blue}[16]}, and a 0-1 law {\color{blue}[26]}. Additional topics of interest are included in {\color{blue}[19,22,23,31,32]}.

In order to construct Kesten's \textit{incipient infinite cluster} for the metric graph GFF, originally seminally established for two-dimensional Bernoulli percolation {\color{blue}[24]}, for the case $d=3$ we study critical percolation clusters of the metric graph GFF by examining continuous time trajectories of the exploration martingale, whose time dependent filtration $\mathcal{I}$ varies with respect to each step of the exploration process. From the initial exploration step, estimates for the quadratic variation, which is expressed as a superposition of exploration martingales, is dependent on the harmonic measure along the boundary of each box explored on  the metric graph associated to $\textbf{Z}^d$, where at each vertex the transition probability is equal to $\frac{1}{2d}$ from that of the simple random walk, and otherwise for the edges behaves like a standard Brownian motion $\widetilde{B}_t$.

Particularly, for behaviors of metric graph GFF level-sets, over appropriate scales, given a finite volume with side length $N$, the harmonic measure, or alternatively, the crossing probability, is dependent upon a sequence of heights $h( N) \longrightarrow 0$. To introduce other relevant properties of the GFF, for a Gaussian process on $G= G(V,E)$ is characterized as having zero mean and covariance between arbitrary points $u$ and $v$, with respective fields $\phi_u$ and $\phi_v$, given by the Green's function of the Laplacian,

\begin{align*}
  \textbf{E} [ \phi_u \phi_v ] = G(u,v) = \textbf{E}_u \big[ \int_0^{\infty} \textbf{1}_{S_t = v} \big]   \text{ , } \\
\end{align*}

\noindent for $d>3$, given a simple random walk $S_t$. For $d=3$, the stopping time of the random walk, instead of being infinite, is given by $\tau  = \mathrm{inf} \{ t : S_t \in \partial V   \} $, which for $d=2$ is a finite box, while for $d \geq 3$, $V = \textbf{Z}^d$. From the GFF on a graph, $\widetilde{B}_t$ yields the GFF on metric graphs (where the metric graph of $G$ is denoted with $\widetilde{G}$) respectively $\widetilde{\phi_u}$ and $\widetilde{\phi_v}$, whose covariance depends on $G(u,v)$ from $G$, is given by,

\begin{align*}
           G(w_1, w_2) =  (1 - r_1)(1-r_2) G(u_1,u_2) + r_1 r_2 G(v_1,v_2) + (1-r_1)r_2 G(u_1,v_2) + \cdots \\ r_1 (1-r_2) G(v_1,u_2) + 2d ( r_1 \wedge r_2 - r_1 r_2) \textbf{1}_{(u_1,v_1) = (u_2,v_2)}    \text{ , } \\
\end{align*}

\noindent where the vertices of the Green’s function are given by adjacent pairs of vertices $(u_1,v_1)$ and $(u_2,v_2)$, with $w_1 \in e(u_1,v_1)$, $w_2 \in e(u_2,v_2)$, $r_1 = |w_1 - u_1|$, $r_2 = |w_2 - u_2|$, and $d \geq 3$ {\color{blue}[25]}. Otherwise, for $d=2$ we enforce the stopping time $\widetilde{\tau} = \mathrm{inf} \{ t :     \widetilde{B_t} \in \partial V   \}$. From the construction of the GFF whose covariance function is respectively defined above for graphs and metric graphs, we introduce previous results on percolation probabilities.

\subsection{This paper's contributions}

The following effort seeks to construct the IIC for the metric graph GFF through estimates on crossing probabilities. Under the assumptions that Kesten initially placed on Bernoulli percolation, with knowledge that the critical point over $\textbf{Z}^2$ equals $\frac{1}{2}$, he demonstrated that one can separate crossings across larger scales in terms of smaller scales; in order to formulate variants of Kesten's seminal arguments for constructing the IIC in two-dimensions, Basu and Sapozhnikov introduce a robust method that can apply not only to higher-dimensional IICs, but also potentially to other two, or three, dimnensional models which have wider possible classes of boundary conditions. Despite the fact that individual bonds in Bernoulli percolation configurations can be closed, or opened, other models in Statistical Mechanics can exhibit significantly different macroscopic properties through constructions of boundary conditions that are far more intricate. In the case of the metric graph GFF, the metric graph extension of the ordinary GFF inherits Dirichlet boundary conditions, which is closely related to the formulation of the GFF in terms of the Dirichlet energy. Moreover, constructing the metric graph extension of the GFF in two, and three, dimensions alike presents more promising opportunities for obtaining crossing probability estimates, in addition to examining other near-critical, and critical, phenomena which can be dependent on the choice of boundary conditions, characteristic length, correlation length, and several related objects. Beyond the forthcoming extension of the IIC construction for the metric graph GFF proposed by Basu and Sapozhnikov, we address a research direction of interest raised in a work of Dubedat and Falconet. That is, as another model of interest, the Villain model exhibits a different encoding of boundary conditions, and of the model itself, rather than those of the ordinary, and metric graph, GFFs. Nevertheless, by introducing a suitable notion of the crossing probability, over the Villain configuration space a summation of exponentials, which will be shown to be related to the transition probability for Brownian motion, enables for a construction of an IIC type limit over the cable system of the Villain model.

\subsection{Statement of previous results on level-set percolation probabilities for the metric graph extension of the GFF} 

We collect necessary results for a quantification of how critical percolation clusters on level-sets from the metric graph GFF. To study connectivity properties in the infinite volume limit, define the box $V_r \equiv [-r,r] \cap \textbf{Z}^d$ centered about the origin for $d \geq 3$ and $r \geq 1$, and is otherwise a finite box for $d=2$. We define the GFF level-sets as a random surface defined over the vertices of $\widetilde{\textbf{Z}^d}$, with $\widetilde{E}^{\geq h}_N = \{   v \in \widetilde{\textbf{Z}^d} : \widetilde{\phi}_{N,v} \geq h    \}$. Finally, we examine the chemical distance between two arbitrary vertices on $\widetilde{Z}^d$, $u$ and $v$, in which $u$ and $v$ are connected with respect to the graph distance if $\mathcal{D}_{N,h} = \mathrm{inf} \{   \mathcal{D}_{N,h}(u,v) : u \in A, v \in B , A , B \subset \widetilde{\textbf{Z}^d}      \} < \infty$, while the vertices are disconnected if $\mathcal{D}_{N,h} = \infty$, indicating that the metric graph GFF does not attain a crossing to the $\partial V$ with height $\geq h$. The percolation probability associated to such a crossing from $u$ to $v$ is introduced below.

\bigskip
\noindent\textbf{Theorem 1} (\textit{percolation probability bounded away from $0$ uniformly in $N$}, {\color{blue}[11, Theorem 1]}) For any height $h \in \textbf{R}$, and $0 < \alpha < \beta < \gamma < 1$, denote $\textbf{P}^{\alpha,\gamma}_{N,h} \big[ \cdot \big]$ as the metric-graph GFF law of $\widetilde{\phi}_N$ conditioned on the crossing event that $u$ and $v$ are connected, namely $\{ V_{\alpha N} \overset{\geq h}{\longleftrightarrow} \partial V_{\gamma N} \} $, with respect to the graph distance $\mathcal{D}_{N,h}$. For any $\epsilon >0$, there exists some $C$ so that in the infinite volume limit, crossings across $V_N$ satisfy,

\begin{align*}
    \underset{{N \longrightarrow +\infty}  }{\mathrm{lim \text{ }  sup}} \text{ }   \textbf{P}^{\alpha, \gamma}_{N,h} [ \mathcal{D}_{N,h}(V_{\alpha N} , \partial V_{\beta N}) > C N (\mathrm{log} N)^{\frac{1}{4}} ] \leq \epsilon     \text{ . } \\
\end{align*}

\noindent\textbf{Theorem 2} (\textit{exponential decay of percolation probabilities; poly-logarithmic upper bound for percolation probabilities above the $h=0$ level-set for three dimensions, and greater than three dimensions}, {\color{blue}[11, Theorems 3 $\&$ 4]}) For a strictly positive sequence of heights $h_N$, satisfying,

\begin{align*}
   h_N \leq C  \text{, } 
\end{align*}

\noindent for sufficient $C>0$, the three dimensional percolation probability is bounded above exponentially, in which there exists some $c$ so that,

\begin{align*}
     p_{N,h} \leq \mathrm{exp}\big(         - \frac{c h^2 N}{\mathrm{log}N}   \big)       \text{ , } \\
\end{align*}

\noindent for $N > 1$, while for all dimensions greater than $3$ the upper bound takes a similar form,

\begin{align*}
      p_{N,h} \leq \mathrm{exp} \big(      -   c h^2 N    \big)     \text{ . } \\
\end{align*}

\noindent Above the $h=0$ level-set, the percolation probabilities to a surface whose boundary is at length $k$ of the exploration process satisfies,

\begin{align*}
 \frac{c}{\sqrt{N}} \leq p_{N,0} \leq C \sqrt{\frac{\mathrm{log N}}{N}}  \text{ , } \\
\end{align*}

\noindent for $d=3$,  and otherwise,

\begin{align*}
      \frac{c}{\sqrt{n}} \leq p_{N,0} \leq \frac{C}{\sqrt{n}}  \text{ , } \\
\end{align*}

\noindent for $d > 3$.

\bigskip

\noindent\textbf{Definition 1} (\textit{exploration martingale}, {\color{blue}[11]}) For vertices $v\in V$ from $\widetilde{G}$, the exploration martingale is the conditional expectation $M_{A,t} = \textbf{E}[X_A | \mathcal{F}_{\mathcal{I}_t}]$, where the observable $X_A$ takes the form,

\begin{align*}
  X_A = \frac{1}{|A|} \sum_{v \in A} \widetilde{\phi}_v  \text{ , }\\
\end{align*}

\noindent where $A \subset V$. The filtration of the conditional expectation satisfies,

\begin{align*}
    \mathcal{F}_{\mathcal{I}_t} =   \{   \mathcal{E} \in \mathcal{F}_{\widetilde{G}} :  \mathcal{E} \cap \{  \mathcal{I}_t \subset U   \} \in \mathcal{F}_U \text{ for all open } \mathcal{I}_0 \subset U          \}      \text{ , } \\
\end{align*}

\noindent where $\mathcal{I}_0$ denotes the initial starting point of the random walk on $h$ before exploration has begun, while $\mathcal{I}_t$ denotes the step of the exploration process of the random walk at $t$. From the observable above, the exploration martingale takes the form,

\begin{align*}
   M_{A,t} = \textbf{E}[X_A | \mathcal{F}_{\mathcal{I}_t}]        \text{ , } \\
\end{align*}

\noindent where in the final quantity in the equality above, $M_{A,t}$ can be decomposed through contributions from a harmonic measure along the boundary of exploration and values of the GFF at vertices of $\partial \mathcal{I}_t$.

\bigskip

\noindent\textbf{Definition 2} (\textit{quadratic variation}, {\color{blue}[11,Corollary 10]}) The quadratic variation of the exploration martingale is defined as,

\begin{align*}
    \langle M_A \rangle_t = \mathrm{Var}[X_A | \mathcal{F}_{\mathcal{I}_0}] - \mathrm{Var}[X_A | \mathcal{F}_{\mathcal{I}_t}]  \text{ . } \\
\end{align*}

\noindent Besides the almost-sure continuity of the exploration martingale and quadratic variation {\color{blue}[11, Lemma 9]}, from the GFF covariance, subtracting the Green's functions associated with random walks at times $s$ and $t$, where $0 \leq s \leq t$, and $v , v^{\prime} \in A$, admits the following harmonic measure decomposition,

\begin{align*}
        G_s(v,v^{\prime}) - G_t(v, v^{\prime}) \equiv    \sum_{w \in \mathcal{I}_t} \mathrm{Hm}_t(v,w) G(w,v^{\prime})   - \sum_{w^{\prime}\in \mathcal{I}_s} \mathrm{Hm}_s(v,w^{\prime}) G(w^{\prime},v^{\prime})   \text{ , } \\
\end{align*}

\noindent from the accompanying sets of explorations on $\widetilde{G}$, $\mathcal{I}_s$ and $\mathcal{I}_t$. Next, we gather results surrounding difference quotients for the Harmonic measure $\mathrm{Hm}$, in addition to the cumulative distribution and survivor functions for the standard normal distribution, each of which we respectively denote with $\Phi$ and $\bar{\Phi}$. 

\bigskip

\noindent\textbf{Definition 3} (\textit{harmonic measure quotients, level-set stopping time}, {\color{blue}[11, Lemma 14]}) Define the two functions,

\begin{align*}
    f_1(N) =  \text{ } \underset{{K \longrightarrow \text{ } + \text{ } \infty}}{\mathrm{lim \text{ } sup}}  \frac{ \mathrm{Hm}(x_K, V_N) - \mathrm{Hm}(x_K, 0) }{G(0,x_K)}   \text{ , }          \\
\end{align*}

\noindent and,

\begin{align*}
         f_2(N) =  \text{ } \underset{{K\longrightarrow \text{ } + \text{ }  \infty}}{\mathrm{lim \text{ } inf}} \inf_{\mathcal{I}}    \frac{ \mathrm{Hm}(x_K, \mathcal{I})-\mathrm{Hm}(x_k,0)}{G(0,x_K)}   \text{ . }       \\
\end{align*}

\noindent Also, for a standard Brownian motion $B$ previously introduced over $G$, the associated stopping time with an arbitrary level-set height $h$ is,

\begin{align*}
       \tau_h = \mathrm{inf} \big\{     \tau \geq 0 : B_t    \leq  ht - \frac{\widetilde{\phi}_0 - h}{\sigma^2_d}           \big\}            \text{ , } \\
\end{align*}

\noindent where at each time step $t$ of evolution $B$ is independent of $\widetilde{\phi}_0$.

\bigskip

\noindent\textbf{Definition 4} (\textit{stopping time-normal distribution equality}, {\color{blue}[11, Proposition 15]}) Fix $m \in \textbf{R}$, $b > 0$, and stopping time $\tau = \mathrm{inf} \{ t>0 : B_t \leq mt - b \}$. For some $T > 0$, the following equality holds,

\begin{align*}
       \textbf{P}[\tau \leq T] = \Phi\big[\frac{b}{\sqrt{t}} - m \sqrt{T}\big] + \mathrm{exp}(2bm) \Phi\big[\frac{b}{\sqrt{T}} + m\sqrt{T}\big]    \text{ . } \\
\end{align*}

\bigskip

\noindent \textbf{Corollary 1} (\textit{stopping time- percolation probability estimates}, {\color{blue}[11, Corollary 16]}) From the quotient of harmonic measures given in \textbf{Definition 3}, the percolation probability to the boundary of a box of length $N$ and level-set $h$ satisfies,

\begin{align*}
    \textbf{P}\big[ \tau_h \geq f_1(N)\big] \leq p_{N,h} \leq \textbf{P}\big[\tau_h \geq f_2(N)\big]   \text{ , } \\
\end{align*}

\noindent from the level-set stopping time.

\subsection{Paper organization}

We study the previously mentioned open problem raised in {\color{blue}[11]} pertaining to the three-dimensional, metric graph GFF, through analysis of the quantities defined in the previous subsection, namely the exploration martingale, and harmonic measure which is introduced as an alternative for quantifying crossing probabilities across level-sets of the metric graph GFF within the critical height window. To do this, we discuss results from seminal arguments due to {\color{blue}[24]} for constructing the IIC for Bernoulli percolation, in which two probabilistic quantities are shown to be equal, one of which is the limit of supercritically defined conditional crossing probabilities, which, through a sequence of parameters $> \frac{1}{2}$ approaching criticality at $p_c \equiv \frac{1}{2}$, is shown to equal another conditional crossing probability taken exactly at $p_c$. For the metric graph GFF, to study criticality about the height threshold about $h \equiv 0$, we make use of arguments due to {\color{blue}[11]} which analyze the behavior of the percolation probability from sequences of levels of the metric graph GFF strictly above, and below, the aforementioned critical height threshold.

\section{Properties of the IIC from Kesten's 1986 work}

 \subsection{Two dimensional Bernoulli percolation objects }

\noindent A streamlined presentation of arguments provided in {\color{blue}[24]} is provided in {\color{blue}[3]}, the key quantities of which are summarized below. From classical assumptions, fix some $p \in [0,1]$, and denote $\textbf{P}_p[ \cdot ]$ as the probability measure for sampling a Bernoulli bond percolation configuration of parameter $p$, with critical parameter,

\begin{align*}
p_c \equiv \mathrm{inf} \big\{ p : \textbf{P}_p \big[ | C( v) | \equiv + \infty      \big] > 0         \big\} \text{, } 
\end{align*}

\noindent which is the smallest parameter for which the cardinality of the open cluster, $C(v)$, for some $v \in V ( \textbf{Z}^2 )$, is infinite with positive probability. Finally, from the induced graph metric $\rho$ for $\textbf{Z}^2$, denote the set,

\begin{align*}
  S ( v, n ) \equiv   \text{ }  \big\{ x \in V (\textbf{Z}^2) :  \rho(v,x) = n   \big\} \text{ , } 
\end{align*}

\noindent for positive integers $m \leq n$.

\bigskip

\noindent \textbf{Definition} \textit{5} (\textit{cylinder events}, {\color{blue}[24]}). A \textit{cylinder event} $E$ is an event that is dependent upon a finite number of vertices from $\textbf{Z}^2$.

\bigskip

\noindent \textbf{Definition} \textit{6} (\textit{conditional crossing probability for Bernoulli percolation at criticality}, {\color{blue}[24]}). Introduce, for $p \equiv p_c$,

\begin{align*}
    \mathcal{P}_1 \equiv   \textbf{P}_{p_c} \big[ E     |    w \longleftrightarrow        S ( w , n)   \big]     \text{. }
\end{align*}

\bigskip

\noindent \textbf{Definition} \textit{7} (\textit{conditional crossing probability when approaching the critical parameter from above}, {\color{blue}[24]}). Introduce, as $p \searrow p_c$, for $w \in V$,

\begin{align*}
  \mathcal{P}^{\prime}_2 \equiv   \underset{p \searrow p_c }{\mathrm{lim}} \mathcal{P}_2 \equiv  \underset{p \searrow p_c }{\mathrm{lim}}  \textbf{P}_p \big[  E   | | C (w) | \equiv + \infty \big]       \text{. }
\end{align*}

\bigskip

\noindent In addition to the \textbf{Definitions} provided above, to demonstrate, beyond the existence of Kesten's IIC about the critical height threshold $h \equiv 0$ of the metric graph GFF, we will also make use of the following assumption.

\bigskip

\noindent \textit{Assumption One} (\textit{Lower bounds for crossing probabilities of Bernoulli percolation at the critical parameter for Kolmogorov extension of the incipient infinite cluster measure to configurations of edges}). At the critical parameter $p_c \equiv \frac{1}{2}$, the following assumption from {\color{blue}[23]}, consisting of two inequalities below, asserts,

\begin{align*}
   \textbf{P}_{p_c} \bigg[          [ - \Lambda , 0 ] \times [0 , n ]      \underset{[-\Lambda , 3n + \Lambda ] \times [0,n] }{\longleftrightarrow}   [3n , 3n + \Lambda ] \times [0,n]   \bigg]    \geq \delta           \text{, } \\     \textbf{P}_{p_c} \bigg[   [0,n ] \times [ - \Lambda , 0 ]          \underset{[0,n] \times [-\Lambda , 3n + \Lambda ]}{\longleftrightarrow}          [0,n] \times [3n , 3n + \Lambda]      \bigg]    \geq \delta               \text{, } 
\end{align*}

\noindent for some sufficient, strictly positive, lower bound $\delta$, and constant $\Lambda \equiv \Lambda (G)$, for which,

\begin{align*}
  \forall e \in G : |e| \leq \Lambda   \text{. }
\end{align*}

\bigskip

\noindent Below, we present a modification of \textit{Assumption One}, with \textit{Assumption Two}. In the following, as will be introduced in \textbf{Definition 9}, denote the metric graph GFF measure, precisely at the critical height threshold, with $\textbf{P}_{h \equiv 0} [ \cdot ] \equiv \textbf{P}_0 [ \cdot ]$, below.

\bigskip

\bigskip

\noindent \textbf{Theorem} \textit{3} (\textit{conditional cylinder events under the Bernoulli bond percolation are equal, under the previously stated Assumption}, {\color{blue}[24]}). Given the assumption provided above, the limits satisfy,

\begin{align*}
    \mathcal{P}_1      \equiv \mathcal{P}^{\prime}_2  \text{. }
\end{align*}

\noindent Denoting the common value shared by the two limits above with $\nu ( E)$, by Kolmogorov extension, there exists only one infinite cluster for which,

\begin{align*}
           \nu \big[\forall w \in V ,  \exists ! \text{ }  C(w) :  | C(w) | = + \infty      \big] \equiv 1            \text{, } 
\end{align*}

\noindent where the unique infinite cluster contains $w$.

\subsection{Three dimensional metric-graph GFF objects}

\noindent From objects introduced in the previous section, we now introduce analogously defined objects for the metric graph GFF. From {\color{blue}[24]} it is known that the critical threshold for the metric graph GFF is $h \equiv 0$. Under this choice of parameters, in three-dimensions critical level-set percolation for the metric graph GFF can be further examined. First, state the following property that the metric graph GFF satisfies, namely the Strong Markov Property.

\bigskip

\noindent \textbf{Theorem} \textit{4} (\textit{GFF Strong Markov Property}, \textbf{Theorem} \textit{8}, {\color{blue}[11]}). For any random compacted connected subset $K \subsetneq \widetilde{G}$, conditionally upon $K$ and the filtration $\mathcal{F}_K$, the following equality in distribution,

\begin{align*}
      \big\{        \widetilde{\phi_v}   : v \in \widetilde{G} \text{ } \backslash \text{ } K       \big\}       \overset{\mathrm{d}}{\equiv}     \big\{   \textbf{E} \big[   \widetilde{\phi_v}                 \text{ }   \big|  \text{ }        \mathcal{F}_K            \big]  +         v \in \widetilde{G} \text{ } \backslash \text{ } K      \big\}    \text{, } 
\end{align*}

\noindent holds, concerning the metric graph GFF $\widetilde{\phi_v}$, for $v \in \widetilde{G}$, over the complementary metric graph $\widetilde{G} \backslash K$ under Dirichlet boundary conditions. From the equality in distribution above, the conditional expectation on the LHS admits the decomposition,

\begin{align*}
    \textbf{E} \big[ \widetilde{\phi_v} | \mathcal{F}_k ]     \equiv   \underset{u \in \partial K}{\sum}  \mathrm{Hm}(v , u ; K) \widetilde{\phi_u}      \text{ } \text{ , }
\end{align*}

\noindent over a summation of contributions from the harmonic measure $\mathrm{Hm}( \cdot , \cdot ; \cdot)$.

\bigskip

\noindent In addition to the property above, it is also of importance to introduce the following quantities from \textbf{Definitions} \textit{8-10}. From the metric graph $\widetilde{G}$, let $\rho$ be the corresponding graph metric, for which, given some $v \in \widetilde{G}$, define,

\begin{align*}
 \mathcal{B} ( v, n ) \equiv    \big\{   x \in V : \rho ( v , x ) \leq n       \big\}            \text{, } 
\end{align*}

\noindent corresponding to the box that is centered about $v$, which is of some strictly positive length $n$, in addition to,

\begin{align*}
 \partial \mathcal{B} (v , n ) \equiv       \big\{      x \in V : \rho ( v , x ) \equiv n         \big\}            \text{, } 
\end{align*}

\noindent corresponding to the boundary of the box introduced in the previous item. Finally, define, for $m \leq n$,

\begin{align*}
    \mathcal{A}(v, m , n ) \equiv   \mathcal{B} ( v , n ) \backslash \mathcal{B} ( v , m-1)  \equiv \big\{ x \in V : m -1 \leq \rho(v,x) \leq n \big\}      \text{, } 
\end{align*}

\noindent corresponding to the finite volume enclosed by annulus contained within $\mathcal{B}(v,n)$, and $\mathcal{B}(v,m)$, respectively.

\bigskip

\noindent \textbf{Definition} \textit{8} (\textit{cylinder events that are a function of the height of the metric graph GFF}). Introduce, for some $h >0$,

\begin{align*}
 \mathscr{E} ( h ) \equiv \mathscr{E}   \text{, } 
\end{align*}

\noindent corresponding to the cylinder event of the metric graph GFF.

\bigskip

\noindent For \textbf{Definition} \textit{9} below, denote another vertex of $\widetilde{G}$ with $w$.

\bigskip

\noindent \textbf{Definition} \textit{9} (\textit{Kesten's IIC for the metric graph GFF in three dimensions, from conditional crossing probabilities for level-set percolation of the metric graph GFF at the critical $h \equiv 0$ threshold; the metric graph GFF analogue to} \textbf{Definition} \textit{6} \textit{originally provided for Bernolli percolation in two dimensions}). Introduce, for $h \equiv 0$,

\begin{align*}
   \textbf{P}^{\alpha, \beta}_{N,0}  \big[   \mathscr{E}  |     w \overset{ \geq 0}{\longleftrightarrow} \mathcal{B}( w , n )       \big]   \equiv  \textbf{P}_{N,0}     \big[  \mathscr{E} |     w \overset{ \geq 0}{\longleftrightarrow} \mathcal{B}( w , n )    \big] \equiv \textbf{P}_0     \big[      \mathscr{E}  |    w \overset{ \geq 0}{\longleftrightarrow} \mathcal{B}( w , n )          \big]  \text{ } \text{ , }
\end{align*}

\noindent corresponding to the level-set percolation crossing event across the critical GFF height threshold, conditionally upon the existence of a configuration such that $\big\{    w \longleftrightarrow \mathcal{B}( w , n )      \big\}$ occurs. 

\bigskip

\noindent For \textbf{Definition} \textit{10} below, denote the open cluster of the metric graph GFF, centered about some $v$, with $\mathcal{C}(v)$.

\bigskip

\noindent \textbf{Definition} \textit{10} (\textit{supercritical conditional crossing probability for level-set percolation connectivity events of the metric graph GFF above the critical} $h \equiv 0$ \textit{threshold, the metric graph GFF analogue to} \textbf{Definition} \textit{7} \textit{originally provided for Bernoulli percolation}). Introduce, for a sequence of heights $\{ h_N \} \longrightarrow 0$, where $N$ is the length of the side of a finite volume, 

\begin{align*}
    \underset{ h_N \searrow 0}{\mathrm{lim}}  \textbf{P}^{\alpha, \beta}_{N,h_N}  \big[   \mathscr{E}  |  | \mathcal{C}(v) | = + \infty       \big]  \equiv    \underset{ h_N \searrow 0}{\mathrm{lim}}  \textbf{P}_{N,h_N}  \big[   \mathscr{E}  |  | \mathcal{C}(v) | = + \infty       \big]   \equiv  \underset{ h_N \searrow 0}{\mathrm{lim}}  \textbf{P}_{h_N}  \big[   \mathscr{E}  |  | \mathcal{C}(v) | = + \infty       \big]   \text{, } 
\end{align*}

\noindent corresponding to the limit of the conditional probability that a cylinder event $\mathscr{E}$ occurs, given the existence of an open cluster, $\mathcal{C}$, centered about $v$, with infinite cardinality. 

\bigskip

\noindent \textbf{Defintion} \textit{11} (\textit{conditional crossing probability for level-set percolation in order to decouple configurations from infinity}). Given a sequence of heights of the metric graph GFF, beginning with $h_0$ from $\big\{ h_N \big\}$, with a corresponding sequence of metric graphs $\{ \widetilde{G}_N \}$, introduce the following quantity, as $h_N \searrow 0$,

\begin{align*}
 \underset{h \in [ 0 , h_0]}{\mathrm{sup}}  \textbf{P}_N \big[        \mathscr{E}_2       |     \mathscr{E}_1                   \big]                  \text{, } 
\end{align*}

\noindent corresponding to the conditional probability of a \textit{cylinder} event $\mathscr{E}_2 \big( h \big) \equiv \mathscr{E}_2$ occurring, conditionally upon the occurrence of $\mathscr{E}_1 \big( h \big) \equiv \mathscr{E}_1$, which are two \textit{cylinder} events, that are respectively given by, 

\begin{align*}
    \mathscr{E}_1 \equiv  \big\{  \mathcal{B}(v , m )     \overset{ \geq h}{\longleftrightarrow }          \mathcal{B} ( v, n )                  \big\}     \text{, } 
\end{align*}

\noindent for $m \leq n$, corresponding to the existence of a crossing event between the boxes $\mathcal{B}(v , m )$ and $ \mathcal{B} ( v, n )$, and by,

\begin{align*}
     \mathscr{E}_2 \equiv    \bigg\{      \forall v \in V , \exists w_1 \neq v , \cdots , w_n \neq w_{n-1} \neq  \cdots \neq w_1   \neq v \in V : \mathcal{C}(v) \cap \bigg( \text{ } \underset{w_1 \leq w \leq w_n }{\bigcap} \text{ }  \mathcal{C}(w) \bigg)  \equiv \emptyset          \bigg\}                   \text{, } 
\end{align*}

\noindent corresponding to the simultaneous existence of at least \textit{two} disjoint open clusters besides the existence of the open cluster $\mathcal{C}(v)$.

\bigskip

\noindent To proceed, in light of the metric graph GFF objects presented thus far, we formulate assumptions that were originally provided for construction of the IIC for Bernoulli percolation in two dimensions, from {\color{blue}[24]}, in terms of the properties that the metric graph GFF would be expected to fulfill, for the construction of the IIC to remain applicable in three dimensions.

\bigskip

\noindent \textbf{Assumptions on the metric graph GFF IIC in three dimensions}. The metric graph GFF IIC is based upon the following two assumptions imposed upon $\widetilde{G}$:

\begin{itemize}
    \item[$\bullet$] \underline{Assumption one} (\textit{uniqueness of metric graph GFF open clusters}): For any height $h$ of the metric graph GFF within the critical window, there exists almost surely at most $1$ infinite open cluster.
    
      \item[$\bullet$] \underline{Assumption two} (\textit{quasi-multiplicativity of crossing probabilities for level-set percolation in the metric graph GFF, which the authors of {\color{blue}[3]} state is expected to hold for low dimensional integer-lattices of dimension $< 6$}): Fix some $v \in V$ and $\delta >0$. There exists a strictly positive constant $\mathscr{C}_{*}$, such that for a range of height parameters of the metric graph free field defined above the critical threshold $h \equiv 0$, given by $ h \in [0, \delta]$, the probability of the level-set connectivity event,
      
      \begin{align*}
     \big\{         \mathscr{X}                \overset{\mathscr{Z}, \geq h}{\longleftrightarrow}              \mathscr{Y}             \big\}      \text{, } 
      \end{align*}

      \noindent occurring admits the following lower bound,

      \begin{align*}
         \textbf{P}_h \big[      \mathscr{X}                \overset{\mathscr{Z}, \geq h}{\longleftrightarrow}              \mathscr{Y}         \big]    \geq              \textbf{P}_h \big[       \mathscr{X}              \overset{\mathscr{Z} , \geq h}{\longleftrightarrow}                   \partial \mathcal{B} (v ,  )        \big]    \mathscr{C}_{*}              \textbf{P}_h \big[                \mathscr{Y}              \overset{\mathscr{Z} , \geq h}{\longleftrightarrow}                   \partial \mathcal{B} (v ,  )                           \big]                      \text{, } 
      \end{align*}

      \noindent for a connected subset $\mathscr{Z} \subset V$, with $\mathscr{Z} \supseteq \mathcal{A}(v, m , )$, $\mathscr{X} \subset \mathscr{Z} \cap \mathcal{B}(v , m)$, $\mathscr{Y} \subset \mathscr{Z} \backslash \mathcal{B} ( v, )$.

\end{itemize}

  \bigskip

  \noindent From the two assumptions provided above, we make use of the FKG property of discrete and continuous Gaussian free fields, as mentioned following the statement of \textbf{Theorem} \textit{1.1} in {\color{blue}[21]}, and also in {\color{blue}[28]}. We denote this property of the GFF as (FKG).

  \bigskip

     \noindent Also, introduce the following \textbf{Theorem}.

     \bigskip

     \noindent \textbf{Theorem} \textit{5} (\textit{establishing equivalence between the random variables provided in} \textbf{Definition} \textit{9}, \textit{and in} \textbf{Definition} \textit{10}, \textit{under appropriate assumptions on the three-dimensional metric graph}). Suppose that $\widetilde{G}$ satisfies \underline{Assumption one} and \underline{Assumption two} provided above, with some choice of $v \in V$ and suitable $\delta > 0$. Then, for any cylinder event $\mathscr{E}$, one obtains the equivalence,

        \begin{align*}
            \textbf{P}_0     \big[      \mathscr{E}  |    w \overset{ \geq 0}{\longleftrightarrow} \mathcal{B}( w , n )          \big]  \equiv       \underset{ h_N \searrow 0}{\mathrm{lim}}  \textbf{P}_{h_N}  \big[   \mathscr{E}  |  | \mathcal{C}(v) | = + \infty       \big]      \text{, } 
        \end{align*}

        \noindent between the conditional crossing event given on the LHS, precisely at $h \equiv 0$, and the conditional crossing event on the RHS, for a sequence of heights $\longrightarrow 0$. Under the equality provided above, by Kolmogorov extension as provided in the second condition of $\textbf{Theorem}$ \textit{3},

\begin{align*}
  \nu_h \big[    \forall w \in \widetilde{V} , \exists ! C ( \widetilde{w})   : | C(\widetilde{w}) | = + \infty       \big]  \equiv 1 \text{, } 
\end{align*}

\noindent for the vertex set $\widetilde{V}$ of the metric graph $\widetilde{G}$.

        \bigskip

        \noindent \textit{Proof of Theorem 5}. We make use of the arguments described initially for Bernoulli percolation, provided in {\color{blue}[24]}, in addition to the modifications to that argument provided in {\color{blue}[3]} for demonstrating that the IIC for Bernoulli percolation also exists on infinite, connected, bounded degree, graphs, with the following. Namely, it suffices to demonstrate, by Kolmogorov extension, that,

        \begin{align*}
        \textbf{P}_h \big[    \mathscr{E} | w \overset{\geq h}{\longleftrightarrow} \mathcal{B}(w,n)      \big]           \longrightarrow \nu_h ( \mathscr{E})   \text{, } \tag{\textit{1}}
        \end{align*}

        \noindent holds for $h \in [0 , \delta]$, and $\delta >0$, as provided in \underline{Assumption two}, in which the above convergence is uniform. Towards accomplishing this, for some increasing event, $i >0$, and parameters satisfying \underline{Assumption two}, denote, for $N_{i+1} \equiv 2^{i+1} \sqrt{\mathrm{log} n_{i+1}} > 2^i \sqrt{\mathrm{log} n_i } \equiv N_i$, $\mathcal{B}_i \equiv \mathcal{B}(v , N_i)$ corresponding to the box centered about $v$ of length $N_i$, $\partial \mathcal{B}_i \equiv \partial \mathcal{B} ( v , N_i)$ corresponding to the boundary of the $\mathcal{B}(v , N_i)$ centered about $v$ of length $N_i$, and $\mathcal{A}_i \equiv \mathcal{A}(v , N_{i+1} , N_i ) \equiv \mathcal{B}(v , N_{i+1}) \backslash \mathcal{B}(v, N_i)$ corresponding to the annulus centered about $v$ of length $N_{i+1} - N_i$. Furthermore, introduce,

        \begin{align*}
             \epsilon_i \equiv   \underset{\mathscr{F}_1 , \mathscr{F}_2 \in \mathcal{A}_i}{\underset{ h \in [0,h_0]}{\mathrm{sup}}}   \textbf{P}_h \big[     \big\{   \mathscr{F}_1 \overset{\geq h}{\underset{\mathcal{A}_i}{\longleftrightarrow}}       \mathscr{F}_2       \big\}^c     |    \partial \mathcal{B}_i \overset{\geq h}{\longleftrightarrow} \partial \mathcal{B}_{i+1}                                   \big]        \text{, }  \tag{\textit{2}}
        \end{align*}
        
        \noindent as the supremum of the conditional crossing probability, as provided above for height parameters taken over the critical window $[0,h_0]$, with $h_0$ denoting the first height in the sequence $\{h_N\}$ which monotonically decreases to the critical threshold $h \equiv 0$. We further manipulate this quantity, observing,
        
           \begin{align*}
                 \textbf{P}_h \big[  \big\{           w  \overset{ \geq h}{\longleftrightarrow} \partial \mathcal{B}  (w,n) \big\} \cap   \big\{   \mathscr{F}_1     \overset{\geq h}{\underset{\mathcal{A}_1}{\longleftrightarrow} }      \mathscr{F}_2 \big\}   \big]   \leq \textbf{P}_h \big[    w \overset{\geq h}{\longleftrightarrow} \partial \mathcal{B}_i                    \big] \text{ } \textbf{P}_h \big[    \big\{ \partial \mathcal{B}_i \overset{\geq h}{\longleftrightarrow} \partial \mathcal{B}_{i+1} \big\} \cap       \big\{   \mathscr{F}_1     \overset{\geq h}{\underset{\mathcal{A}_1}{\longleftrightarrow}}       \mathscr{F}_2 \big\}   \big] \text{ } \times \cdots \\ \textbf{P}_h \big[      \partial \mathcal{B}_{i+1} \overset{\geq h}{\longleftrightarrow} \partial \mathcal{B}(w,n)         \big]             \\ \leq     \epsilon_i           \textbf{P}_h \big[   w \overset{\geq h}{\longleftrightarrow} \partial \mathcal{B}_i    \big]    \textbf{P}_h \big[   \partial \mathcal{B}_i \overset{\geq h}{\longleftrightarrow} \partial \mathcal{B}_{i+1}           \big] \textbf{P}_h \big[   \partial \mathcal{B}_{i+1} \overset{\geq h}{\longleftrightarrow} \partial \mathcal{B}(w,n)    \big]   \\ \leq  \epsilon_i  \big( \mathscr{C}_{*} \big)^2     \textbf{P}_h \big[ w \overset{\geq h}{\longleftrightarrow} \partial \mathcal{B} ( w, n ) \big]       \text{, }  \tag{$1^{*}$}
        \end{align*}
        
        \noindent where the square of $\mathscr{C}_{*}$, in the uppermost bound, appears in the statement of \underline{Assumption two} for quasi-multiplicativity. Pursuant of the decomposition steps provided for constructing the IICs for Bernoulli percolation that are first provided in {\color{blue}[24]}, and subsequently in {\color{blue}[3]}, introduce,
        
        \begin{align*}
          \mathcal{C}_i \equiv \big\{   x \in \mathcal{B}_i : x \overset{\geq h}{\underset{\mathcal{B}_{i+1}}{\longleftrightarrow}} \mathcal{B}_i         \big\}           \text{, } 
        \end{align*}

       \noindent corresponding to points within $\mathcal{B}_{i+1}$ that are connected to $\mathcal{B}_i$, in addition to, from the denomination $\partial \mathcal{B}_{i+1} \equiv \partial \mathcal{B} ( v, N_i + 1)$,

       \begin{align*}
         \mathcal{D}_i \equiv \big\{  \forall   x \in \partial{B}_{i+1}  ,  \exists x \sim  y \in \mathcal{C}_i     : \big\{x \overset{  \ngeq h }{\longleftrightarrow} y\big\}    \big\}   \text{, } 
       \end{align*}
       
       \noindent corresponding to the points in $\partial{B}_{i+1}$ for which there exists an edge $(xy)$ which does not contain a crossing of at least height $h_N$. From these two random objects, introduce, for $\mathscr{U} \subset \mathcal{B}_i$ and $\mathscr{R} \subset \partial \mathcal{B}_{i+1}$,
       
       \begin{align*}
     \mathcal{F}_i \big( \mathscr{U}  ,  \mathscr{R}  \big) \equiv \big\{    \mathcal{C}_i \equiv \mathscr{U}      \big\} \cup \big\{     \mathcal{D}_i \equiv \mathscr{R}     \big\}   \text{, } 
       \end{align*}

       \noindent corresponding to the event that $\mathcal{C}_i \equiv \mathscr{U}$, and that $\mathcal{D}_i \equiv \mathscr{R}$ simultaneously. From all admissible pairs $(\mathscr{U} , \mathscr{R})$, given $\mathcal{C}_i$ and $\mathscr{F}_i (\mathscr{U} , \mathscr{R} \big) \neq \emptyset$ for which $\big\{ \mathcal{C}_i \equiv \mathscr{U} \big\}$, introduce the union,
       
       \begin{align*}
  \mathcal{F}_i \big( \mathscr{U} , \mathscr{V} , \underset{i}{\prod} \big) \equiv \mathcal{F}_i  \equiv      \underset{(\mathscr{U} , \mathscr{R}) \in \prod_i}{\bigcup}               \mathcal{F}_i \big( \mathscr{U} ,   \mathscr{R}   \big)           \text{, } 
       \end{align*}
       
       \noindent over the space $\prod_i$  of all admissible $( \mathscr{U} , \mathscr{R} )$. For suitable $n$ with $n > N_{i+1} + N_0$,
       
       \begin{align*}
           \textbf{P}_h \big[ \mathscr{E} \cap  \big\{           w  \overset{\geq h}{\longleftrightarrow} \partial \mathcal{B}  (w,n) \big\} \cap   \big\{   \mathscr{F}_1     \underset{\mathcal{A}_1}{\overset{\geq h}{\longleftrightarrow}}       \mathscr{F}_2 \big\}   \big] \equiv    \underset{(\mathscr{U} , \mathscr{R}) \in \prod_i}{\sum} \textbf{P}_h \big[    \mathscr{E} \cap \big\{ w \overset{\geq h}{\longleftrightarrow} \partial \mathcal{B}(w, n ) \big\} \cap \big\{ \big\{    \mathcal{C}_i \equiv \mathscr{U}       \big\} \cup \big\{     \mathcal{D}_i \equiv \mathscr{R}     \big\}    \big\}    \big]                     \text{, } 
       \end{align*}
       
       \noindent which, upon making the observation that,
       
       \begin{align*}
     \textbf{P}_h \big[    \mathscr{E} \cap \big\{ w \overset{\geq h}{\longleftrightarrow} \partial \mathcal{B}(w, n ) \big\} \cap \big\{ \big\{    \mathcal{C}_i \equiv \mathscr{U}       \big\} \cup \big\{     \mathcal{D}_i \equiv \mathscr{R}     \big\}    \big\}    \big] \equiv \textbf{P}_h \big[      \mathscr{E} \cap \big\{  w \overset{\geq h}{\longleftrightarrow} \partial \mathcal{B}_{i+1}    \big\} \cap                \big\{ \big\{    \mathcal{C}_i \equiv \mathscr{U}       \big\} \cup \big\{     \mathcal{D}_i \equiv \mathscr{R}     \big\}    \big\}    \big]   \times \cdots \\ \textbf{P}_h \big[        \mathscr{R} \overset{\geq h}{\underset{\mathcal{B} ( w, n )  \backslash {\mathscr{U} }}{\longleftrightarrow}}    \partial \mathcal{B}_n   \big]   \text{, } 
       \end{align*}
       
        \noindent is equivalent to the following summation over $( \mathscr{U} , \mathscr{V} )$,
        
         \begin{align*}
      \underset{(\mathscr{U} , \mathscr{R}) \in \prod_i}{\sum}   \textbf{P}_h \big[      \mathscr{E} \cap \big\{  w \overset{\geq h}{\longleftrightarrow} \partial \mathcal{B}_{i+1}    \big\} \cap                \big\{ \big\{    \mathcal{C}_i \equiv \mathscr{U}       \big\} \cup \big\{     \mathcal{D}_i \equiv \mathscr{R}     \big\}    \big\}    \big]    \textbf{P}_h \big[        \mathscr{R} \overset{\geq h}{\underset{\mathcal{B} ( w, n )  \backslash \mathscr{U} }{\longleftrightarrow}}    \partial \mathcal{B}_n   \big]   \text{. }
       \end{align*}
       
       \noindent To quantify the absolute value between random variables that have been previously manipulated, form the difference, between the probability of the intersection of the increasing event $\mathscr{E}$ occurring, intersected with,
       
       \begin{align*}
        \mathscr{E} \cap \big\{w \overset{\geq h}{\longleftrightarrow}  \partial \mathcal{B}(w, n) \big\} \text{, } 
       \end{align*}

       \noindent and each possible realization of,

       \begin{align*}
       \mathscr{E} \cap \big\{ w \overset{\geq h}{\longleftrightarrow}     \partial \mathcal{B}_{i+1} \big\} \cap \big\{ \big\{ \mathcal{C}_i \equiv \mathscr{U} \big\} \cup \big\{ \mathcal{D}_i \equiv \mathscr{R} \big\} \big\}           \text{, } 
       \end{align*}

       \noindent  for admissible $(\mathscr{U} , \mathscr{V}) \in \prod_i$ occurring, with, 
       
       \begin{align*}
     \big|   \textbf{P}_h \big[ \mathscr{E} \cap \big\{ w \overset{\geq h}{\longleftrightarrow} \partial \mathcal{B}(w, n ) \big\} \big]                 -  \underset{(\mathscr{U} , \mathscr{V}) \in \prod_i}{\sum}      \textbf{P}_h \big[      \mathscr{E} \cap \big\{  w \overset{ \geq h}{\longleftrightarrow} \partial \mathcal{B}_{i+1}    \big\} \cap               \underset{\mathcal{F}_i}{\underbrace{\big\{ \big\{    \mathcal{C}_i \equiv \mathscr{U}       \big\} \cup \big\{     \mathcal{D}_i \equiv \mathscr{R}     \big\}    \big\} } }  \big]  \times \cdots \\   \textbf{P}_h \big[        \mathscr{R} \overset{ \geq h}{\underset{\mathcal{B} ( w, n )  \backslash \mathscr{U} }{\longleftrightarrow}}    \partial \mathcal{B}_n   \big]            \big|  \tag{\textit{DIFF}} \\ \leq  \text{ }  \textbf{P}_h \big[ \mathscr{E} \cap   \mathcal{F}_i     \cap                \big\{ w   \overset{\geq h}{\longleftrightarrow}  \partial \mathcal{B}(w ,n) \big\}   \big]     \\  \overset{ (
     \mathrm{QM})}{\leq}  \textbf{P}_h \big[  \mathcal{F}^c_i  \cap \big\{   w \overset{ \geq h }{\longleftrightarrow} \partial \mathcal{B}(w,n) \big\}  \big] \big(\mathscr{C}_{*}\big)^{-2} \textbf{P}_h \big[  \mathscr{E}  \big] \text{ }  \text{ , }  \tag{$2^{*}$}
     \end{align*}
     
     \noindent where the final expression, upon making the observation that the following conditional probability,
     
     \begin{align*}
{\textbf{P}_h \big[ \mathcal{F}^c_i | w \overset{\geq h}{\longleftrightarrow} \partial \mathcal{B}(w,n)] }\equiv \frac{\textbf{P}_h \big[    \mathcal{F}^c_i \cap \big\{ w \overset{\geq h}{\longleftrightarrow} \partial \mathcal{B}(w,n) \big\}     \big]}{\textbf{P}_h \big[    \mathcal{F}^c_i    \big]}         \text{, }    \tag{$3^{*}$}
     \end{align*}

     \noindent is equivalent to the ratio of probabilities given on the RHS above, yields the expression,
     
     \begin{align*}
\textbf{P}_h \big[    \mathcal{F}^c_i    \big]   \textbf{P}_h \big[    \mathcal{F}^c_i |  w \overset{\geq h}{\longleftrightarrow} \partial \mathcal{B}(w,n)     \big]  \text{ } \big(\mathscr{C}_{*}\big)^{-2}  \text{ }  \textbf{P}_h \big[     \mathscr{E}   \big] \\ \underset{\mathscr{F}_1 , \mathscr{F}_2 \in \mathcal{A}_i}{\equiv}  \textbf{P}_h \big[    \mathcal{F}^c_i    \big]   \text{ } \textbf{P}_h \big[   \big\{  \mathscr{F}_1  \overset{\geq h}{\underset{\mathcal{A}_i}{\longleftrightarrow}} \mathscr{F}_2        \big\}^c      | w \overset{\geq h}{\longleftrightarrow}    \partial \mathcal{B}(w,n)        \big]  \big(\mathscr{C}_{*}\big)^{-2} \textbf{P}_h \big[ \mathscr{E} \big]   \\  \overset{(*)}{\leq} \epsilon^{\prime}_i \big(\mathscr{C}_{*}\big)^{-2}  \textbf{P}_h \big[ \mathscr{E} \big]  \\  \equiv \epsilon^{\prime}_i \big( \mathscr{C}_{*} \big)^{-2} \frac{\textbf{P}_h \big[ \mathscr{E} \big] \textbf{P}_h \big[  w \longleftrightarrow \partial \mathcal{B}(w,n) \big] }{\textbf{P}_h \big[  w \longleftrightarrow \partial \mathcal{B}(w,n) \big]}  \\  \overset{(\mathrm{FKG})}{\leq} \epsilon^{\prime}_i \big(\mathscr{C}_{*}\big)^{-2} \frac{\textbf{P}_h \big[ \mathscr{E} \cap \big\{  w \overset{\geq h}{\longleftrightarrow}  \partial \mathcal{B}(w,n) \big\} \big]}{\textbf{P}_h \big[     w \longleftrightarrow \partial \mathcal{B}(w,n)  \big] } \\ \overset{(**)}{\leq} \epsilon^{\prime} \big( \mathscr{C}^{\prime}_{*} \big)^{-2} \textbf{P}_h \big[ \mathscr{E} \cap \big\{  w \overset{\geq h}{\longleftrightarrow}  \partial \mathcal{B}(w,n) \big\} \big] \text{, } 
       \end{align*}
       
        \noindent where in $(*)$, the constant in the upper bound satisfies,

        \begin{align*}
       \epsilon^{\prime}_i  \geq \epsilon_i \text{ } \textbf{P}_h \big[ \mathcal{F}^c_i \big]     \text{, } 
        \end{align*}

        \noindent while in $(**)$,

        \begin{align*}
            \big(\mathscr{C}^{\prime}_{*}\big)^{-2}   \geq \frac{\big( \mathscr{C}_{*} \big)^{-2}}{\textbf{P}_h \big[  w \longleftrightarrow \partial \mathcal{B}(w,n)  \big]}   \text{, } 
        \end{align*}

        \noindent in addition to the fact that the cylinder event $\mathscr{E}$, and,
        
        \begin{align*}
         \big\{ w \overset{\geq h}{\longleftrightarrow}   \partial \mathcal{B}(w,n) \big\}       \text{, } 
        \end{align*}

        \noindent being increasing events implicates the following lower bound,
        
        \begin{align*}
        \textbf{P}_h \big[  \mathscr{E} \cap \big\{ w \overset{ \geq h}{\longleftrightarrow}  \partial \mathcal{B}(w,n) \big\}   \big]   \geq     \textbf{P}_h \big[  \mathscr{E}  \big] \textbf{P}_h \big[      w \overset{ \geq h}{\longleftrightarrow}  \partial \mathcal{B}(w,n)         \big]    \text{. }
        \end{align*}
        
        \noindent

        \noindent To proceed from final upper bound obtained from (\textit{**}), introduce,

        \begin{align*}
    u^{\prime} ( \mathscr{U} , \mathscr{R} , h ) \equiv       u_h^{\prime} ( \mathscr{U} , \mathscr{R}  )    \equiv  \textbf{P}_h \big[  \mathscr{E} \cap \big\{   w    \overset{ \geq h}{\longleftrightarrow} \partial \mathcal{B}_{i+1} \big\} \cap \mathcal{F}_i  \big]                      \text{, }  \\       u^{\prime\prime} ( \mathscr{U} , \mathscr{R} , h  ) \equiv        u_h^{\prime\prime} ( \mathscr{U} , \mathscr{R} )  \equiv  \textbf{P}_h \big[ \big\{   w \overset{\geq h}{\longleftrightarrow}  \partial \mathcal{B}_{i+1} \big\} \cap \mathcal{F}_i  \big]     \text{, } \\    \gamma ( \mathscr{U} , \mathscr{R} , n , h) \equiv   \gamma_h ( \mathscr{U} , \mathscr{R} , n)  \equiv \textbf{P}_h \big[    \mathscr{R} \underset{\mathcal{B}(w,n) \backslash \mathscr{U}}{\overset{ \geq h} {\longleftrightarrow} }  \partial \mathcal{B}(w,n)\big]   \text{, } 
        \end{align*}

        \noindent which respectively corresponding to the probability of the intersection of the three events provided occurring, the probability of the intersection of the two events occurring, and finally, the probability of a connection between $\mathscr{R} \longleftrightarrow \partial \mathcal{B}(w,n)$ occurring, within $\mathcal{B}(w,n) \backslash \mathscr{U}$. Fix $(\mathscr{U} , \mathscr{R}) \in \prod_i$, $j > i+2$, and $n > N_{j+1} + N_0$. Observe that one can form a similar expression to $(\textit{DIFF})$ with,

        \begin{align*}
     \big|   \gamma ( \mathscr{U} , \mathscr{R} , n , h )        -    \underset{(\mathscr{U}^{\prime} , \mathscr{R}^{\prime}) \in \prod_j }{\sum}      \textbf{P}_h \big[    \big\{   \mathscr{R} \underset{\mathcal{B}_{j+1} \backslash \mathscr{U}}{\overset{\geq h}{\longleftrightarrow}}   \partial \mathcal{B}_{j+1} \big\} \cap \mathcal{F}_{j-1} \cap \mathcal{F}_j \big( \mathscr{U}^{\prime} , \mathscr{R}^{\prime}  \big)   \big]               \big|  \text{, } 
        \end{align*}
        
      \noindent which, upon similiar rearrangements as given previously, yields the upper bound,
      
      \begin{align*}
   \big( \mathscr{C}^{\prime}_{*}\big)^{-2}  \bigg( \underset{\epsilon_{j-1}}{\underbrace{\underset{\mathscr{F}_1 , \mathscr{F}_2 \in \mathcal{A}_i}{\underset{ h \in [0,h_0]}{\mathrm{sup}}}   \textbf{P}_h \big[     \big\{   \mathscr{F}_1 \overset{\geq h}{\underset{\mathcal{A}_i}{\longleftrightarrow}}       \mathscr{F}_2       \big\}^c     |    \partial \mathcal{B}_{j-1} \overset{\geq h}{\longleftrightarrow} \partial \mathcal{B}_{j}                                   \big] }}   +  \underset{\epsilon_{j}}{\underbrace{\underset{\mathscr{F}_1 , \mathscr{F}_2 \in \mathcal{A}_i}{\underset{ h \in [0,h_0]}{\mathrm{sup}}}   \textbf{P}_h \big[     \big\{   \mathscr{F}_1 \overset{\geq h}{\underset{\mathcal{A}_i}{\longleftrightarrow}}       \mathscr{F}_2       \big\}^c     |    \partial \mathcal{B}_{j+1} \overset{\geq h}{\longleftrightarrow} \partial \mathcal{B}_{j+2}                                   \big]  }}  \bigg) \times \cdots \\      \textbf{P}_h \big[  \mathscr{R} \underset{\mathcal{B}(w,n) \backslash \mathscr{U}}{\overset{\geq h}{\longleftrightarrow} } \partial \mathcal{B}(w,n)  \big]              \text{. }
      \end{align*}
      
      \noindent Furthermore, executing this procedure repeatedly, fix some $k,N>0$, with $k, N \in \textbf{N}$. For a sequence of parameters $j_1 , \cdots , j_N$, with $j_{k+1} > j_k +2$, and as under a previous assumption, $n > N_{j_k+1} + N_0$, that,
      
      \begin{align*}
       \underset{(\mathscr{U}^{\prime} , \mathscr{V}^{\prime}) \in \prod_j}{\sum}              \gamma ( \mathscr{U} , \mathscr{R} , n , h )  \text{ } \text{ , }
      \end{align*}
      
      \noindent for the conditional probability $\mathscr{M}$ appearing in the summation above, which is defined with,
      
      \begin{align*}
    \mathscr{M} \big( \mathscr{U} , \mathscr{R} ; \mathscr{U}^{\prime} , \mathscr{R}^{\prime} \big)  \equiv \textbf{P}_h \big[   \big\{      \mathscr{R} \underset{\mathcal{B}_{j+1} \backslash \mathscr{U}}{\overset{\geq h}{\longleftrightarrow}}   \partial \mathcal{B}_{i+1} \big\} \cap     \mathcal{F}_{j-1} \cap \mathcal{F}_j \big( \mathscr{U}^{\prime} , \mathscr{R}^{\prime}  \big)  \big]       \text{, }     \tag{$4^{*}$}
      \end{align*}
      
      \noindent is bounded from below with,
      
      \begin{align*}
   \gamma( \mathscr{U} , \mathscr{R} , n , h ) -   \big( \mathscr{C}^{\prime}_{*} \big)^{-2} \big( \epsilon_{j-1} + \epsilon_j \big)   \gamma( \mathscr{U} , \mathscr{R} , n , h )                             \equiv   \gamma( \mathscr{U} , \mathscr{R} , n , h ) \big( 1 - \big( \mathscr{C}^{\prime}_{*} \big)^{-2} \big( \epsilon_{j-1} + \epsilon_j \big) \big)              \text{, } 
      \end{align*}
      
      \noindent and above with,
      
      \begin{align*}
        \gamma( \mathscr{U} , \mathscr{R} , n , h ) +  \big( \mathscr{C}^{\prime}_{*} \big)^{-2} \big( \epsilon_{j-1} + \epsilon_j \big)   \gamma( \mathscr{U} , \mathscr{R} , n , h )                             \equiv   \gamma( \mathscr{U} , \mathscr{R} , n , h ) \big( 1 + \big( \mathscr{C}^{\prime}_{*} \big)^{-2} \big( \epsilon_{j-1} + \epsilon_j \big) \big)                   \text{, } 
      \end{align*}
      
      \noindent hence constituting a tight bound on,

     \begin{align*}
        \underset{(\mathscr{U}^{\prime} , \mathscr{V}^{\prime}) \in \prod_j}{\sum}              \gamma ( \mathscr{U} , \mathscr{R} , n , h )   \text{. }
     \end{align*}

      \noindent Next, observe, for $j_{k+1}$ introduced previously, the sequence of pairs, from the collection,
      
      \begin{align*}
    \big\{    \big( \mathscr{U}^{\prime}_{j_1} , \mathscr{R}^{\prime}_{j_1} \big) \in \prod_{j_1} , \cdots , \big(   \mathscr{U}^{\prime}_{j_{k+1}} , \mathscr{R}^{\prime}_{j_{k+1}} \big)  \in \prod_{j_{k+1}}   \big\} \text{, } 
      \end{align*}
      
      \noindent which implies, for $h \in [ 0 , h_0]$,

      \begin{align*}
    \frac{u_h^{\prime} ( \mathscr{U} , \mathscr{R} , h)}{u_h^{\prime\prime} ( \mathscr{U} , \mathscr{R} , h) }   \underset{(   \mathscr{U}^{\prime}_{j_{k+1}} , \mathscr{R}^{\prime}_{j_{k+1}} )  \in \prod_{j_{k+1}}  }{\underset{\vdots}{\underset{ ( \mathscr{U}^{\prime}_{j_1} , \mathscr{R}^{\prime}_{j_1} ) \in \prod_{j_1} }{\underset{2\leq k\leq N}{\prod}} }}  \frac{\mathscr{M}_h \big( \mathscr{U}_1 , \mathscr{R}_1 ; \mathscr{U}^{\prime}_2 , \mathscr{R}^{\prime}_2 \big) }{\gamma_h \big( \mathscr{U}_k , \mathscr{R}_k ; \mathscr{U}^{\prime}_k , \mathscr{R}^{\prime}_k \big) }                    \frac{\mathscr{M}_h\big(   \big(\mathscr{U}_1\big)_k , \big(\mathscr{R}_1\big)_k    ; \big( \mathscr{U}^{\prime}_2\big)_k , \big( \mathscr{R}^{\prime}_2\big)_k  \big) }{\mathscr{M}_h \big(   \big(\mathscr{U}_1\big)_{k-1} , \big(\mathscr{R}_1\big)_{k-1}    ; \big( \mathscr{U}^{\prime}_2\big)_{k-1} , \big( \mathscr{R}^{\prime}_2\big)_{k-1}    \big) }                     \times \cdots \\  \big( \frac{\mathscr{M}_h \big( \mathscr{U}_1 , \mathscr{R}_1 ; \mathscr{U}^{\prime}_2 , \mathscr{R}^{\prime}_2 \big)}{\gamma_h\big( \mathscr{U}_k , \mathscr{R}_k ; \mathscr{U}^{\prime}_k , \mathscr{R}^{\prime}_k \big) } \big)^{-1}             \text{, } 
      \end{align*}

      \noindent is tightly bound from below with,

      \begin{align*}
         \mathrm{exp} \big( -\epsilon \big) \textbf{P}_h \big[ \mathscr{E} | w \overset{ \geq h}{\longleftrightarrow} \partial \mathcal{B} (w, n ) \big]          \text{, } 
      \end{align*}

      \noindent and is tightly bound from above with,
      
      \begin{align*}
           \mathrm{exp} \big( \epsilon \big) \textbf{P}_h \big[ \mathscr{E} | w \overset{ \geq h}{\longleftrightarrow} \partial \mathcal{B} (w, n ) \big]                  \text{, } 
      \end{align*}

      \noindent for suitable $\epsilon > 0$. As demonstrated in {\color{blue}[24]}, and subsequently adapted for establishing the existence of the IIC for Bernoulli percolation in {\color{blue}[3]} under different assumptions on the underlying graph structure, it suffices to argue than an inequality of the form,
      
      \begin{align*}
        \big( \mathscr{M}_h \big(    \mathscr{U}_1 , \mathscr{R}_1 ; \mathscr{U}^{\prime}_1 , \mathscr{R}^{\prime}_1    \big)     \text{ } \mathscr{M}_h \big(    \mathscr{U}_2 , \mathscr{R}_2 ; \mathscr{U}^{\prime}_2 , \mathscr{R}^{\prime}_2    \big)     \big) \text{ } \big(  \mathscr{M}_h   \big(    \mathscr{U}_1 , \mathscr{R}_1 ; \mathscr{U}^{\prime}_2 , \mathscr{R}^{\prime}_2                            \big)   \mathscr{M}_h   \big(  \mathscr{U}_2 , \mathscr{R}_2 ; \mathscr{U}^{\prime}_1 , \mathscr{R}^{\prime}_1                         \big)   \big)^{-1}                   \leq \kappa^2       \text{, }  
      \end{align*}

      \noindent holds, in which the upper bound is independent of $\epsilon$ and $k$. To demonstrate that the upper bound provided above holds, consider similar random sets that were introduced previously, each of which are given by,

      \begin{align*}
       \mathcal{X}_j   =   \big\{        x \in \mathcal{A}_{j-1} : x \underset{\mathcal{A}_{j-1}}{\overset{\geq h}{{\longleftrightarrow}}}  \partial \mathcal{B}_j  \big\} \text{, }  \\     \mathcal{Y}_j   =   \big\{     \forall y \in \partial \mathcal{B}(v , N_{j-1} -1 )                    \text{ }   , \text{ }     \exists x \sim  y    : \big\{x \overset{  \ngeq h }{\longleftrightarrow} y\big\}    \big\}    \text{, } 
      \end{align*}

      \noindent which can be used to formulate the event, as similarly provided in the union of two events for defining $\mathcal{F}_i$, with, 
      
      \begin{align*}
     \mathcal{G}_j \big(   X                , Y  \big) \equiv             \big\{   \mathcal{X}_j \equiv X   \big\} \cup \big\{     \mathcal{Y}_j \equiv Y   \big\}    \text{, } 
      \end{align*}

      \noindent for suitable $X$ and $Y$. From this choice of $X$ and $Y$, observe that the occurrence of $\mathcal{G}_j$ provided above is dependent upon the edges within $\mathcal{A}_{j-1}$, in which it is possible that,
      
      \begin{align*}
            \big\{ \mathcal{X}_j \equiv X \big\} \cap \mathcal{F}_{j-1} \equiv \emptyset     \text{, } 
      \end{align*}
      
      \noindent or that,
      
        \begin{align*}
               \big\{ \mathcal{X}_j \equiv X \big\} \subset  \mathcal{F}_{j-1}        \text{. }
      \end{align*}

      \noindent From the objects defined above, for each admissible pair $(X,Y)$, $(\mathscr{U}, \mathscr{R}) \in \prod_i$, $(\mathscr{U}^{\prime} , \mathscr{R}^{\prime} ) \in \prod_j$, $\mathcal{F}_{j-1}$ is expressible with,

      \begin{align*}
  \mathcal{F}_{j-1} \equiv          \underset{(X, Y) \in \Gamma_j}{\bigcup}    \mathcal{G}_j \big( X , Y \big)           \text{, } 
      \end{align*}

      \noindent for $\Gamma_j$, which is expressible with another union,
      
      \begin{align*}
         \Gamma_j \equiv \underset{\mathcal{G}_j ( X , Y) \neq \emptyset}{\underset{\{ \mathcal{X}_j \equiv X \} \subset \mathcal{F}_{j-1} }{\bigcup}} \big\{ X , Y \big\}             \text{, } 
      \end{align*}
      
      \noindent which can then be used to define the following summation which is a product over the probabilities,
      
      \begin{align*}
         \mathscr{M} \big( \mathscr{U} , \mathscr{R} ; \mathscr{U}^{\prime} , \mathscr{R}^{\prime}  , h \big)    \equiv     \mathscr{M}_h \big( \mathscr{U} , \mathscr{R} ; \mathscr{U}^{\prime} , \mathscr{R}^{\prime}  \big)   \equiv \underset{(X,Y) \in \Gamma_j}{\sum}            \textbf{P}_h \big[   \mathscr{R} \underset{\mathcal{B}_j \backslash ( X \cup \mathscr{U})}{\overset{\geq h}{\longleftrightarrow}} Y          \big] \times \cdots \\    \textbf{P}_h \big[           \mathcal{G}_j ( X, Y)   \cap \mathcal{F}_j \big( \mathscr{U}^{\prime} , \mathscr{R}^{\prime} \big) \cap \big\{   Y     \overset{\geq h}{\longleftrightarrow}  \mathscr{R}^{\prime}    \big\}                                   \big]  \\ \overset{(\mathrm{FKG})}{\geq}   \underset{(X,Y) \in \Gamma_j}{\sum}            \textbf{P}_h \big[   \mathscr{R} \underset{\mathcal{B}_j \backslash ( X \cup \mathscr{U})}{\overset{\geq h}{\longleftrightarrow}} Y          \big]    \textbf{P}_h \big[           \mathcal{G}_j ( X, Y) \big] \textbf{P}_h \big[    \mathcal{F}_j \big( \mathscr{U}^{\prime} , \mathscr{R}^{\prime} \big) \big] \times \cdots \\ \textbf{P}_h \big[   Y     \overset{\geq h}{\longleftrightarrow}  \mathscr{R}^{\prime}                                    \big]              \text{. }
      \end{align*}

      \noindent From the summation of probabilities above, the first term,

      \begin{align*}
      \textbf{P}_h \big[ \mathscr{R} \underset{\mathcal{B}_j \backslash ( X \cup \mathscr{U} ) }{\overset{\geq h}{\longleftrightarrow}} Y      \big]     \text{, }  \tag{$5^{*}$}
      \end{align*}
      
      \noindent can be bound from below, by quasi-multiplicativity of crossing probabilities,

       \begin{align*}
      \textbf{P}_h \big[ \mathscr{R} \underset{\mathcal{B}_j \backslash ( X \cup \mathscr{U})}{\overset{\geq h}{\longleftrightarrow}} Y      \big]  \geq \textbf{P}_h \big[    \mathscr{R} \underset{\mathcal{B}(v , 2 N_{i+1}) \backslash \mathscr{U}} {\overset{\geq h}{\longleftrightarrow}}          \partial \mathcal{B}(v , 2 N_{i+1} )          \big]    \text{ } \mathscr{C}_{*} \text{ }         \textbf{P}_h \big[          \mathcal{G}_j ( X , Y ) \cap \mathcal{F}_j ( \mathscr{U}^{\prime} , \mathscr{R}^{\prime} ) \cap \big\{       Y               {\overset{\geq h}{\longleftrightarrow}} \mathscr{R}^{\prime}   \big\}          \big]  \text{, }      \tag{$6^{*}$}    
      \end{align*}

      \noindent in which the final probability obtained from the quasi-multiplicativity lower bound can be further bound from below with,
      
      \begin{align*}
          \mathscr{C}_{*} \textbf{P}_h \big[   \mathscr{R}        \underset{\mathcal{B}_j \backslash (X \cup \mathscr{U} )}{\overset{\geq h}{\longleftrightarrow}}        Y       \big]              \text{. }
      \end{align*}

    \noindent To obtain the desired lower and upper bounds for concluding the argument, observe, for $\xi \leq 1$,
    
    \begin{align*}
     \big|   \xi    \big|    \big|     \frac{u^{\prime}}{u^{\prime\prime}}         \big| \big|  \xi^{-1}        \underset{2 \leq j \leq k }{\prod}         \frac{\mathscr{M}_j             \big( \mathscr{U}_j \big)        }{\mathscr{M}_j           \big( \mathscr{U}_j \big)        }                  - \frac{u^{\prime\prime}}{u^{\prime}}  \big| \equiv \big|     \frac{u^{\prime}}{u^{\prime\prime}}     \underset{2 \leq j \leq k}{\prod}      \frac{\mathscr{M}_j             \big( \mathscr{U}_j \big)        }{\mathscr{M}_j           \big( \mathscr{U}_j \big)        }  - \xi      \big|    \leq \big( \frac{\kappa-1}{\kappa+1 } \big)^{k-1}         \text{, } 
    \end{align*}

    \noindent for the row vector,

    \begin{align*}
      u^{\prime} \equiv             ( u^{\prime}(1) , \cdots , u^{\prime} (\lambda) )                      \text{, } 
    \end{align*}

    \noindent for the row vector,
    
    \begin{align*}
      u^{\prime\prime} \equiv      ( u^{\prime\prime}(1) , \cdots , u^{\prime\prime} (\lambda) )               \text{, } 
    \end{align*}

        \noindent and for finite volumes $\mathscr{U}_j \subset \mathcal{A}(j_k)$, $\forall j_k$, with $\lambda> 0$. This implies that the lower bound takes the form,
        
          \begin{align*}
  \mathrm{exp} \big( - \epsilon)   \xi -  \mathrm{exp} \big( - \epsilon)     \big( \frac{\kappa-1}{\kappa+1 } \big)^{k-1}      \equiv    \mathrm{exp} \big( - \epsilon)   \big(        \xi -  \big( \frac{\kappa-1}{\kappa+1 } \big)^{k-1}   \big)    \text{, } 
      \end{align*}

        \noindent while the upper bound takes the form,
        
          \begin{align*}
    \mathrm{exp} \big( \epsilon )  \xi + \mathrm{exp} \big( \epsilon ) \big( \frac{\kappa-1}{\kappa+1 } \big)^{k-1}  \equiv  \mathrm{exp} \big( \epsilon )   \big(  \xi +  \big( \frac{\kappa-1}{\kappa+1 } \big)^{k-1}     \big)     \text{, } 
      \end{align*}

        \noindent implying that a bound for,
        
        \begin{align*}
         \textbf{P}_h \big[        \mathscr{E} | w \longleftrightarrow \partial \mathcal{B}(w,n) \big]    \text{, } 
        \end{align*}
        
        \noindent holds, as,

          \begin{align*}
        \mathrm{exp} \big( - \epsilon)   \big(        \xi -  \big( \frac{\kappa-1}{\kappa+1 } \big)^{k-1}   \big)  \leq    \textbf{P}_h \big[         \mathscr{E} | w \longleftrightarrow \partial \mathcal{B}(w,n)      \big]       \leq  \mathrm{exp} \big(  \epsilon)   \big(        \xi + \big( \frac{\kappa-1}{\kappa+1 } \big)^{k-1}   \big)         \text{. }
      \end{align*}
        
        \noindent Hence, for $m , n > N_{i+1} + N_0$, the absolute value difference,      
        
        \begin{align*}
    \big|    \frac{\textbf{P}_h \big[   \mathscr{E}  |  w \overset{\geq h}{\longleftrightarrow} \partial \mathcal{B}(w,m)   \big]  }{\textbf{P}_h \big[   \mathscr{E}   |  w \overset{\geq h}{\longleftrightarrow} \partial \mathcal{B}(w,n) \big]  }       - 1 \big|   \equiv \big|           \textbf{P}_h \big[   \mathscr{E}  |  w \overset{\geq h}{\longleftrightarrow} \partial \mathcal{B}(w,m)   \big]      -   \textbf{P}_h \big[   \mathscr{E}  |  w \overset{\geq h}{\longleftrightarrow} \partial \mathcal{B}(w,n)   \big]      \big| \text{, } 
        \end{align*}
      
      \noindent can be upper bounded with,
      
      \begin{align*}
             \big( \mathrm{exp} \big( \epsilon \big)   - \mathrm{exp} \big( - \epsilon \big)     \big) + \big(  \mathrm{exp} \big(   \epsilon   \big) +  \mathrm{exp} \big(    - \epsilon  \big)     \big)   \big( \frac{\kappa-1}{\kappa+1 } \big)^{k-1}               \text{, } \tag{$7^{*}$}
      \end{align*}
      
      \noindent which implies the desired item in the statement of \textbf{Theorem} \textit{5}, as,
      
      \begin{align*}
      \textbf{P}_0     \big[      \mathscr{E}  |    w \overset{ \geq 0}{\longleftrightarrow} \mathcal{B}( w , n )          \big]  \equiv       \underset{ h_N \searrow 0}{\mathrm{lim}}  \textbf{P}_{h_N}  \big[   \mathscr{E}  |  | \mathcal{C}(v) | = + \infty       \big]    \Rightarrow   \nu_h \big[    \forall w \in \widetilde{V} , \exists ! C ( \widetilde{w})   : | C(\widetilde{w}) | = + \infty       \big]  \equiv 1     \text{, } 
      \end{align*}

      \noindent from which we conclude the proof. \boxed{}
      
      \section{Quasi-mulitiplicativity of crossing probabilities for level-set, metric graph GFF percolation }
      
      \noindent In \textit{Section 3}, with existence of the IIC established in \textit{Setion 2}, we apply arguments developed in {\color{blue}[3]} for proving that quasi-multiplicativity of crossing probabilities for level-set percolation of the metric graph GFF, as for two dimensional Bernoulli percolation, also hold. First, we begin by introducing the setting over which arguments for establishing quasi-multiplicativity take place.
      
      \bigskip

      \noindent \textbf{Definition} \textit{12} (\textit{the three-dimensional environment}). Denote $\mathcal{S} \equiv \big( \mathcal{S} \big( \mathcal{V}\big) , \mathcal{S} \big( \mathcal{E} \big) \big) \equiv \big( \mathcal{V} , \mathcal{E} \big) \subsetneq \textbf{Z}^3$, with,
      
      \begin{align*}
     \mathcal{S} \equiv     \textbf{Z}         \times \big\{       0 , \cdots , k       \big\}^{d-2}     \text{, } 
      \end{align*}
      
      \noindent for $k \geq 0$, and $d \geq 2$.

      \bigskip

     \noindent Besides three-dimensional objects introduced above, we also introduce similar notions as introduced previously for the statement of quasi-multiplicativity, namely a box over $\mathcal{S}$, the boundary of a box over $\mathcal{S}$, and the annulus over $\mathcal{S}$.
     
     \bigskip
     
     \noindent \textbf{Definition} \textit{13} (\textit{three-dimensional boxes, boundaries of three-dimensional boxes, and three-dimensional annuli, for further analysis of crossing events}). Denote the box, with strictly positive side length $2n$, over $\textbf{Z}^3$, with,
     
     \begin{align*}
  \mathcal{B}_{\mathcal{S}}    \big( v , n  \big)  \equiv      \big\{ x \in \mathcal{V} \big( [-n,n]^2 \big)  : \rho_{\mathcal{S}} \big( v , x) \leq n \big\}   \text{, } 
     \end{align*}

     \noindent for the restriction of the graph $\rho_{\mathcal{S}}$, with $\rho$ on $\mathcal{V} \subsetneq \mathcal{V}$. As a result, denote the boundary with,

     \begin{align*}
      \partial \mathcal{B}_{\mathcal{S}} \big( v , n \big) \equiv   \partial \mathcal{B}_{\mathcal{S}} \big( n \big) \equiv  \big\{   x \in     \mathcal{V} \big( [-n,n]^2 \big)   : \rho_{\mathcal{S}} (v,x)= n    \big\}    \text{, } 
     \end{align*}
     
     \noindent and the \textit{three-dimensional annulus} with, for $m > n > 0$,
     
     \begin{align*}
      \mathcal{A}_{\mathcal{S}} (v , m , n ) \equiv  \mathcal{A}_{\mathcal{S}} ( m , n ) \equiv       \mathcal{B}_{\mathcal{S}}(v , n ) \backslash \mathcal{B}_{\mathcal{S}}(v, m -1 )         \equiv \big\{  x \in \mathcal{V} \big( [-n,n]^2 \big) : m-1 \leq \rho_{\mathcal{S}} \big( v , x  \big)  \leq n         \big\}   \text{. }
     \end{align*}
      
   \bigskip
   
   \noindent With the geometric objects defined above over $\mathcal{S}$, also introduce the analogue of quasi-multiplicativity of crossing probabilities over slabs.
   
   \bigskip
      
      \noindent \textbf{Theorem} \textit{6} (\textit{establishing equivalence between the metric-graph GFF limits from } \textbf{Theorem} \textit{5 in the three-dimensional environment}). Fix $\mathcal{S}$, $d \geq 2$, $ k \geq 0$, and $\delta < 1$. Then, from the two limits introduced in \textbf{Theorem} \textit{5}, namely,

      \begin{align*}
        \textbf{P}_0     \big[      \mathscr{E}  |    w \overset{ \geq 0}{\longleftrightarrow} \mathcal{B}( w , n )          \big]                \text{, }  
     \end{align*}

      \noindent corresponding to the metric-graph GFF probability measure at the critical $h$ parameter, and,

        \begin{align*}
            \underset{ h_N \searrow 0}{\mathrm{lim}}  \textbf{P}_{h_N}  \big[   \mathscr{E}  |  | \mathcal{C}(v) | = + \infty       \big]           \text{, }  
      \end{align*}
      
      \noindent corresponding to the limit of a sequence subcritically dependent metric-graph GFF probability measures from a sequence of heights, $\{h_N \}$, each exist and are equal.

      \bigskip
      
      \noindent To demonstrate that the result above holds for the metric-graph IIC, introduce the following.
  
     \bigskip
     
     \noindent \textbf{Lemma} \textit{1} (\textit{three-dimensional quasi-multiplicativity}). From the choice of $\mathcal{S}$, $d$, $k$, and $\delta$ introduced in \textbf{Theorem} \textit{6}, additionally introduce heights $ h \in [ 0 , h_0]$, $m \in \textbf{Z}_{+}$, and a strictly positive $\mathcal{C}_{*}$, from which,

      \begin{align*}
      \textbf{P}_h \big[            \mathscr{X} {\overset{\mathscr{Z} , \geq h}{\longleftrightarrow}}   \mathscr{Y} \big]       \geq     \textbf{P}_h \big[         \mathscr{X} \overset{  \mathscr{Z}  , \geq h}{\longleftrightarrow} \partial \mathcal{B}_{\mathcal{S}} \big(   2 \sqrt{\mathrm{log}m}       \big)     \big]  \text{ } \mathcal{C}_{*} \text{ }   \textbf{P}_h \big[             \mathscr{Y} \overset{  \mathscr{Z}    , \geq h}{\longleftrightarrow}   \partial \mathcal{B}_{\mathcal{S}}    \big(    2 \sqrt{\mathrm{log}m}       \big)   \big]    \text{, } 
      \end{align*}

    \noindent holds, given $\mathscr{Z} \subset \mathcal{S}$, $\mathscr{Z} \supset \mathcal{A} ( m , 2 \sqrt{\mathrm{log} m } + 1 )$, $\mathscr{X} \subsetneq \mathscr{Z} \cap \partial \mathcal{B}_{\mathcal{S}} \big( m \big)$, and $\mathscr{Y} \subset  \mathscr{Z} \backslash \mathcal{B}_{\mathcal{S}} \big(          2 \sqrt{\mathrm{log} m } + 1    \big) $.

        \bigskip
        
      \noindent \textit{Proof of Lemma 1}. To establish quasi-multiplicativity in three-dimensions, instead of directly establishing the inequality between crossing probabilities above, observe, that it suffices to demonstrate,

      \begin{align*}
           \textbf{P}_h \big[            \mathscr{X} {\overset{\mathscr{Z} , \geq h}{\longleftrightarrow}}   \mathscr{Y} \big]       \geq     \textbf{P}_h \big[         \mathscr{X} \overset{   \mathscr{Z} , \geq h}{\longleftrightarrow} \partial \mathcal{B}_{\mathcal{S}} \big(     2 \sqrt{\mathrm{log}m} + 1     \big)     \big]  \text{ } \mathcal{C}_{*} \text{ }   \textbf{P}_h \big[             \mathscr{Y} \overset{   \mathscr{Z} , \geq h}{\longleftrightarrow}   \partial \mathcal{B}_{\mathcal{S}}    \big(     2 \sqrt{\mathrm{log}m}      \big)   \big]          \text{, } 
      \end{align*}
      
      \noindent in which the first probability due to the application of quasi-multiplicativitiy above, 
      
      \begin{align*}
          \textbf{P}_h \big[         \mathscr{X} \overset{  \mathscr{Z}  , \geq h}{\longleftrightarrow} \partial \mathcal{B}_{\mathcal{S}} \big(     2 \sqrt{\mathrm{log}m} + 1     \big)     \big]    \text{, } 
      \end{align*}

      \noindent can itself be lower bound by another application of quasi-multiplicativity, as,

      \begin{align*}
            \textbf{P}_h \big[         \mathscr{X} \overset{   \mathscr{Z} , \geq h}{\longleftrightarrow} \partial \mathcal{B}_{\mathcal{S}} \big(     2 \sqrt{\mathrm{log}m} + 1     \big)     \big]  \geq     \textbf{P}_h \big[    \mathscr{X}                   \overset{   \mathscr{Z} , \geq h}{\longleftrightarrow}       \partial \mathcal{B}_{\mathcal{S}} \big( 2 \sqrt{\mathrm{log}m}      \big] \text{ } \mathscr{C}_{*} \text{ }    \textbf{P}_h \big[       \partial \mathcal{B}_{\mathcal{S}} \big(  \frac{4}{3} \sqrt{\mathrm{log}m}       \big)          \overset{   \mathscr{Z} , \geq h}{\longleftrightarrow}        \partial \mathcal{B}_{\mathcal{S}} \big(    2 \sqrt{\mathrm{log}m}  + 1   \big)      \big]    \text{. }
      \end{align*}
      
      \noindent From the second application of quasi-multiplicativity above, adapting results from {\color{blue}[4]}, yields, for $\textbf{P}_h \big[ \cdot \big]$,
      
      \begin{align*}
      \textbf{P}_h \big[    \partial \mathcal{B}_{\mathcal{S}} \big(   \frac{4}{3} \sqrt{\mathrm{log}m}  \big) \overset{\mathscr{Z} , \geq h}{\longleftrightarrow}  \partial \mathcal{B}_{\mathcal{S}} \big(     2 \sqrt{\mathrm{log} m}     \big)       \big] \geq  \textbf{P}_0 \big[   \partial \mathcal{B}_{\mathcal{S}} \big(   \frac{4}{3} \sqrt{\mathrm{log}m}  \big) \overset{\mathscr{Z}, \geq h}{\longleftrightarrow}  \partial \mathcal{B}_{\mathcal{S}} \big(     2 \sqrt{\mathrm{log} m}     \big)         \big]     \text{, } 
      \end{align*}
      
       \noindent for $h \in [0,h_0]$. Hence, from adaptations of arguments provided in {\color{blue}[4]} for the \textit{three-dimensional environment}, one has,

       \begin{align*}
          \textbf{P}_0 \big[   \partial \mathcal{B}_{\mathcal{S}} \big(   \frac{4}{3} \sqrt{\mathrm{log}m}  \big) \overset{\mathscr{Z}, \geq h}{\longleftrightarrow}  \partial \mathcal{B}_{\mathcal{S}} \big(     2 \sqrt{\mathrm{log} m}     \big)         \big]   > 0              \text{. }
       \end{align*}
       
       \noindent Next, consider,

       \begin{align*}
        \textbf{P}_h \big[  \big\{  \mathscr{X} \overset{ \mathscr{Z},\geq h}{\longleftrightarrow}     \partial \mathcal{B}_{\mathcal{S}} \big( 2 \sqrt{\mathrm{log}m} + 1 \big) \big\} \cap \big\{     \mathscr{Y} \overset{\mathscr{Z},\geq h}{\longleftrightarrow}    \partial \mathcal{B}_{\mathcal{S}} \big( 2 \sqrt{\mathrm{log}m}  ) \big\} \cap \mathcal{E}_1   \big]    \text{, }  \tag{$\mathcal{E}_1$ \textit{probability}}
       \end{align*}
        
       \noindent where the event $\mathcal{E}_1$, is given by,

       \begin{align*}
         \mathcal{E}_1 \equiv \bigg\{               \mathcal{V}_1       \underset{\mathcal{A}_{\mathcal{S}}( 2 \sqrt{\mathrm{log}m}   ,  2 \sqrt{\mathrm{log}m} + 1 ) }{\overset{< h}{\longleftrightarrow}}  \mathcal{V}_2  \bigg\}           \text{, } 
       \end{align*}

       \noindent where the vertices in the connectivity event defined above, within the annulus $\mathcal{A}_{\mathcal{S}}( 2 \sqrt{\mathrm{log}m}   ,  2 \sqrt{\mathrm{log}m} + 1 )$, are drawn from the vertex set,

       \begin{align*}
       \mathcal{V}_1 , \mathcal{V}_2 \in \mathcal{V} \big( \mathcal{S} \big)      \text{, } 
       \end{align*}
       
       \noindent the collection of all vertices from the \textit{three-dimensional environment}. This intersection of crossing events, admits, by (FKG), the lower bound,

         \begin{align*}
        \textbf{P}_h \big[       \mathscr{X} \overset{ \mathscr{Z},\geq h}{\longleftrightarrow}     \partial \mathcal{B}_{\mathcal{S}} \big( 2 \sqrt{\mathrm{log}m} + 1 \big)                \big]          \textbf{P}_h \big[        \mathscr{Y} \overset{\mathscr{Z},\geq h}{\longleftrightarrow}    \partial \mathcal{B}_{\mathcal{S}} \big( 2 \sqrt{\mathrm{log}m}  )        \big]             \textbf{P}_h \big[   \mathcal{E}_1   \big]      \text{. }
       \end{align*}

\bigskip

\noindent From the lower bound above, furthermore, observe, that the intersection of events given with the probability below,
       
       \begin{align*}
       \textbf{P}_h \big[     \big\{  \mathscr{X}     \overset{\mathscr{Z} , \geq h }{\longleftrightarrow}    \bar{\Gamma}     \big\} \cap \big\{      \mathscr{Y}  \overset{\mathscr{Z} , \geq h }{\longleftrightarrow}   \bar{\Gamma}    \big\} \cap \mathcal{E}_1                              \big]     \text{, }  \tag{$\bar{\Gamma}$ \textit{probability}}
       \end{align*}

       \noindent admits a bound from below, of the form, by ($\mathrm{FKG}$),

       \begin{align*}
              \textbf{P}_h \big[     \mathscr{X}     \overset{\mathscr{Z} , \geq h }{\longleftrightarrow}    \bar{\Gamma}     \big]        \textbf{P}_h \big[  \mathscr{Y}  \overset{\mathscr{Z} , \geq h }{\longleftrightarrow}   \bar{\Gamma}    \big]   \textbf{P}_h \big[  \mathcal{E}_1          \big]          \text{, } 
       \end{align*}
       
       \noindent for a path $\bar{\Gamma} \neq \Gamma$, satisfying $\bar{\Gamma} > \Gamma$, where $\Gamma$ denotes the minimal open circuit, given by the path infimum,
       
       \begin{align*}
         \Gamma \equiv \underset{v_0 , \cdots , v_n \in \mathcal{V}}{\underset{\gamma \equiv \{  v_0  \overset{< h}{\longleftrightarrow}  v_n \} }{\underset{\mathrm{circuits} \text{ } \gamma }{\mathrm{inf}}} } \big\{   \rho_{\mathcal{S}}     \big( \gamma ,     \partial \mathcal{B}_{\mathcal{S}} \big( 2 \sqrt{\mathrm{log}m} \big)         \big)        \big\}          \text{, } 
       \end{align*}

       \noindent which, as an infimum over circuits $\gamma$ restricted within the annulus $\mathcal{A}_{\mathcal{S}} ( 2 \sqrt{\mathrm{log}m},2 \sqrt{\mathrm{log}m}+1 ) $, for which,
       
       \begin{align*}
         \gamma \cap \mathcal{A}_{\mathcal{S}} ( 2 \sqrt{\mathrm{log}m},2 \sqrt{\mathrm{log}m}+1 ) \neq \emptyset  \text{, } 
       \end{align*}
       
       \noindent are nonempty subsets over $\mathcal{S}$, as of admissible $\gamma$ satisfy,
       
       \begin{align*}
   \gamma \equiv \underset{i: 0 \leq i \leq n }{\bigcup} \gamma_i \equiv \underset{i: 0 \leq i \leq n }{\bigcup} \big\{  \big\{ v_{i} \overset{< h}{\longleftrightarrow}  v_{i+1} \big\} \big\}   \text{, } 
       \end{align*}
       
       \noindent $\forall i$. The deterministic ordering of circuits is obtained from the relation,
       
       \begin{align*}
        \Gamma_1 \geq \Gamma_2 \Longleftrightarrow   \big\{  \forall v_1 , v_2 \in \mathcal{V} \big( \mathcal{S} \big) , \exists v_1 \in \Gamma_1 , v_2 \in \Gamma_2 :    \rho_{\mathcal{S}}    \big(   v_1          , \partial  \mathcal{B}_{\mathcal{S}} \big( 2 \sqrt{\mathrm{log}m} \big)      \big)    \geq   \rho_{\mathcal{S}}   \big(    v_2   ,   \partial  \mathcal{B}_{\mathcal{S}} \big( 2 \sqrt{\mathrm{log}m} \big)      \big)     \big\} \text{, } 
       \end{align*}

       \noindent for two circuits $\Gamma_1$ and $\Gamma_2$, in which for two paths contained around $\mathcal{B}$. For a subset $W \neq \emptyset$, furthermore, denote,
       
       \begin{align*}
\mathcal{S}  \supset W   \supsetneq \bar{W}  \equiv \big\{    z \equiv \big( z_1 , \cdots , z_d \big) \in \mathcal{S} : \big( z_1 , z_2 , x_3 , \cdots , x_d \big) \in W \text{ for some } x_3 , \cdots , x_d                              \big\}    \text{, } 
       \end{align*}
       
       \noindent corresponding to a nonempty subset of the \textit{three-dimensional environment} for which there exists points $x_3 , \cdots , x_d \in W$. From quantities introduced above relating to the ordering of paths in $\mathcal{S}$, express a lower bound for ($\bar{\Gamma}$ \textit{probability}), in which, from previous manipulation,

       \begin{align*}
        (\bar{\Gamma}\text{ }  \textit{probability}) \geq     \textbf{P}_h \big[       \mathscr{X} \overset{ \mathscr{Z},\geq h}{\longleftrightarrow}     \partial \mathcal{B}_{\mathcal{S}} \big( 2 \sqrt{\mathrm{log}m} + 1 \big)                \big]          \textbf{P}_h \big[        \mathscr{Y} \overset{\mathscr{Z},\geq h}{\longleftrightarrow}    \partial \mathcal{B}_{\mathcal{S}} \big( 2 \sqrt{\mathrm{log}m}  )        \big]             \textbf{P}_h \big[   \mathcal{E}_1   \big]   \overset{(\mathrm{FKG})}{\leq} ( \mathcal{E}_1\text{ }  \textit{probability})       \text{. }
       \end{align*}
       
      \noindent As a result, it suffices to demonstrate, for strictly positive $\mathcal{C}_1$,

      \begin{align*}
       \textbf{P}_h \big[          \big\{  \mathscr{X}           \overset{\mathscr{Z},\geq h}{\longleftrightarrow}   \bar{\Gamma} \big\} \cap \big\{   \mathscr{Y}             \overset{\mathscr{Z},\geq h}{\longleftrightarrow}           \bar{\Gamma}       \big\} \cap  \big\{           \mathscr{X}   \overset{ \mathscr{Z} , \geq h}{\not\longleftrightarrow}     \mathscr{Y}        \big\} \cap \mathcal{E}_1        \big] \geq  \mathcal{C}_1 \textbf{P}_h \big[      \mathscr{X}        \overset{\mathscr{Z},\geq h}{\longleftrightarrow}   \mathscr{Y}                \big]   \text{, } \tag{$\bar{\Gamma}$ \text{ } \textit{Probability I}}
      \end{align*}
      
      \noindent in which the lower bound probability consists of solely the pushforward, under $\textbf{P}_h \big[ \cdot \big]$, of,
      
      \begin{align*}
        \big\{   \mathscr{X} \overset{\mathscr{Z} , \geq h}{\longleftrightarrow}     \mathscr{Y}  \big\}       \text{ }\text{ , } 
      \end{align*}

      \noindent in comparison to the upper bound,  which consists of the intersection, pushed forwards under $\textbf{P}_h \big[ \cdot \big]$, of,
      
      \begin{align*}
       \big\{  \mathscr{X}           \overset{\mathscr{Z},\geq h}{\longleftrightarrow}   \bar{\Gamma} \big\} \cap \big\{   \mathscr{Y}             \overset{\mathscr{Z},\geq h}{\longleftrightarrow}           \bar{\Gamma}       \big\} \cap \mathcal{E}_1         \text{ }\text{ , } 
      \end{align*}
      
       \noindent with the intersection of the three events provided above. To rewrite ($\bar{\Gamma}$ \textit{Probability I}), observe,

       \begin{align*}
   \big\{     \mathscr{X} \overset{\mathscr{Z} , \geq h}{\not\longleftrightarrow}   \mathscr{Y}   \big\}     \cap \big\{    \mathscr{Y} \overset{\mathscr{Z} , \geq h}{\not\longleftrightarrow}   \Gamma    \big\}         \subsetneq    \big\{             \mathscr{X}      \overset{\mathscr{Z} , \geq h}{\longleftrightarrow}  \bar{\Gamma }   \big\}   \text{, } \\    \big\{  \mathscr{X}     \overset{\mathscr{Z} , \geq h}{\not\longleftrightarrow}  \Gamma  \big\}       \cap   \big\{  \mathscr{Y} \overset{\mathscr{Z} , \geq h}{\longleftrightarrow}   \Gamma      \big\}           \subsetneq    \big\{             \mathscr{Y}      \overset{\mathscr{Z} , \geq h}{\longleftrightarrow}      \bar{\Gamma}    \big\}    \text{, } \\   \big\{ \mathscr{X}  \overset{\mathscr{Z} , \geq h}{\longleftrightarrow}       \Gamma \big\}   \cap   \big\{ \mathscr{Y}  \overset{\mathscr{Z} , \geq h}{\not\longleftrightarrow}            \Gamma   \big\}              \subsetneq    \big\{    \mathscr{X}                     \overset{\mathscr{Z} , \geq h}{\not\longleftrightarrow}   \mathscr{Y}           \big\}  \cap \mathcal{E}_1  \text{, } 
       \end{align*}
    
    \noindent in which, from the three consecutively defined events provided above, passing to the pushforward under the measure $\textbf{P}_{h} \big[ \cdot \big]$ yields,

    \begin{align*}
      \textbf{P}_h \big[                 \big\{     \mathscr{X} \overset{\mathscr{Z} , \geq h}{\not\longleftrightarrow}   \mathscr{Y}   \big\}     \cap \big\{    \mathscr{Y} \overset{\mathscr{Z} , \geq h}{\not\longleftrightarrow}   \Gamma    \big\}         \big]  \leq  \textbf{P}_h \big[       \mathscr{X}      \overset{\mathscr{Z} , \geq h}{\longleftrightarrow}  \bar{\Gamma }      \big] \text{, } \\    \textbf{P}_h \big[         \big\{  \mathscr{X}     \overset{\mathscr{Z} , \geq h}{\not\longleftrightarrow}  \Gamma  \big\}       \cap   \big\{  \mathscr{Y} \overset{\mathscr{Z} , \geq h}{\longleftrightarrow}   \Gamma      \big\}    \big] \leq  \textbf{P}_h \big[      \mathscr{Y}      \overset{\mathscr{Z} , \geq h}{\longleftrightarrow}      \bar{\Gamma}         \big] \text{, } \\  \textbf{P}_h \big[   \big\{ \mathscr{X}  \overset{\mathscr{Z} , \geq h}{\longleftrightarrow}       \Gamma \big\}   \cap   \big\{ \mathscr{Y}  \overset{\mathscr{Z} , \geq h}{\not\longleftrightarrow}            \Gamma   \big\}        \big] \leq  \textbf{P}_h \big[            \big\{    \mathscr{X}                     \overset{\mathscr{Z} , \geq h}{\not\longleftrightarrow}   \mathscr{Y}           \big\}  \cap \mathcal{E}_1    \big] \text{, } 
    \end{align*}

    \noindent which, upon further inspection, implies,

          \begin{align*}
        \textbf{P}_h \big[ \mathscr{X} \overset{\mathscr{Z} , \geq h}{\not\longleftrightarrow}   \mathscr{Y}  \big] \text{ } \textbf{P}_h \big[     \mathscr{Y} \overset{\mathscr{Z} , \geq h}{\longleftrightarrow}     \Gamma   \big] \overset{(\mathrm{FKG})}{\leq}  \textbf{P}_h \big[                 \big\{     \mathscr{X} \overset{\mathscr{Z} , \geq h}{\not\longleftrightarrow}   \mathscr{Y}   \big\}     \cap \big\{    \mathscr{Y} \overset{\mathscr{Z} , \geq h}{\not\longleftrightarrow}   \Gamma    \big\}         \big]    \text{, }    \\   \textbf{P}_h \big[  \mathscr{X} \overset{\mathscr{Z} , \geq h}{\not\longleftrightarrow}   \Gamma       \big]       \textbf{P}_h \big[  \mathscr{Y} \overset{\mathscr{Z} , \geq h}{\longleftrightarrow}   \Gamma    \big]       \overset{(\mathrm{FKG})}{\leq}                         \textbf{P}_h \big[         \big\{  \mathscr{X}     \overset{\mathscr{Z} , \geq h}{\not\longleftrightarrow}  \Gamma  \big\}       \cap   \big\{  \mathscr{Y} \overset{\mathscr{Z} , \geq h}{\longleftrightarrow}   \Gamma      \big\}    \big]                               \text{, }  \\  \textbf{P}_h \big[   \mathscr{X}                     \overset{\mathscr{Z} , \geq h}{\longleftrightarrow}   \Gamma      \big]  \textbf{P}_h \big[  \mathscr{Y} \overset{\mathscr{Z} , \geq h}{\longleftrightarrow}  \Gamma     \big]   \overset{(\mathrm{FKG})}{\leq}   \textbf{P}_h \big[            \big\{    \mathscr{X}                     \overset{\mathscr{Z} , \geq h}{\longleftrightarrow}   \Gamma          \big\}  \cap \big\{ \mathscr{Y} \overset{\mathscr{Z} , \geq h}{\longleftrightarrow}  \Gamma   \big\}   \big]  \text{, } 
    \end{align*}

    \noindent in which the LHS lower bound of each inequality can be further bound from below by applying $(\mathrm{FKG})$ three times. As a result of these facts, introduce another closely related to ($\bar{\Gamma}$ \textit{Probability I}), for strictly positive, suitable $\mathcal{C}_1$, as previously provided, in which,
    
    \begin{align*}
        \textbf{P}_h \big[ \mathcal{E}_1 \cap \mathscr{E}_1  \cap \mathscr{E}_2 \cap \mathscr{E}_3     \cap \mathscr{E}_4 \cap \mathscr{E}_5                        \big] \leq                           \mathcal{C}_1            \textbf{P}_h \big[ \mathscr{X} \overset{\mathscr{Z} , \geq h}{\longleftrightarrow}   \mathscr{Y}          \big]                                   \text{, } \tag{$\bar{\Gamma}$ \textit{Probability II}}
    \end{align*}

    \noindent for the following five events,
    
    \begin{align*}
      \mathscr{E}_1 \equiv \big\{   \mathscr{X}  \overset{\mathscr{Z} ,\geq h}{\longleftrightarrow}    \bar{\Gamma}           \big\}   \text{, } \\ \mathscr{E}_2 \equiv \big\{ \mathscr{Y} \overset{\mathscr{Z} , \geq h}{\longleftrightarrow}  \bar{\Gamma}  \big\}   \text{, } \\  \mathscr{E}_3 \equiv  \big\{   \mathscr{X}     \overset{\mathscr{Z} , \geq h}{\not\longleftrightarrow}   \mathscr{Y} \big\}  \\      \mathscr{E}_4 \equiv  \big\{  \mathscr{X} \overset{\mathscr{Z} , \geq h}{\not\longleftrightarrow}  \Gamma \big\}           \text{, } \\   \mathscr{E}_5 \equiv    \big\{       \mathscr{Y} \overset{\mathscr{Z} , \geq h}{\not\longleftrightarrow} \Gamma   \big\}        \text{, } 
    \end{align*}
    
       \noindent in which, in comparison to the intersection of three crossing events,
       
       \begin{align*}
                   \big\{  \mathscr{X}           \overset{\mathscr{Z},\geq h}{\longleftrightarrow}   \bar{\Gamma} \big\} \cap \big\{   \mathscr{Y}             \overset{\mathscr{Z},\geq h}{\longleftrightarrow}           \bar{\Gamma}       \big\} \cap  \big\{           \mathscr{X}   \overset{ \mathscr{Z} , \geq h}{\not\longleftrightarrow}     \mathscr{Y}        \big\}        \text{ } \text{ , }
       \end{align*}

       \noindent that is intersected with $\mathcal{E}_1$ in ($\bar{\Gamma}$ \textit{Probability I}), satisfies,
       
       \begin{align*}
                            \textbf{P}_h \big[ \mathcal{E}_1 \big]   \textbf{P}_h \big[             \mathscr{X} \overset{\mathscr{Z} , \geq h}{\longleftrightarrow}  \mathscr{Y}  \big]  \textbf{P}_h   \big[     \mathscr{Y} \overset{\mathscr{Z} , \geq h}{\longleftrightarrow}  \bar{\Gamma}     \big] \textbf{P}_h \big[  \mathscr{X} \overset{\mathscr{Z} , \geq h}{\not\longleftrightarrow}      \mathscr{Y}     \big] \textbf{P}_h \big[  \mathscr{X} \overset{\mathscr{Z} , \geq h}{\not\longleftrightarrow}   \Gamma   \big] \textbf{P}_h \big[  \mathscr{Y} \overset{\mathscr{Z} , \geq h}{\not\longleftrightarrow}        \Gamma    \big]                            \\  \overset{(\mathrm{FKG})}{\leq }\textbf{P}_h \big[ \mathcal{E}_1 \big]   \textbf{P}_h \big[             \big\{ \mathscr{X} \overset{\mathscr{Z} , \geq h}{\longleftrightarrow}  \mathscr{Y}  \big\} \cap \big\{     \mathscr{Y} \overset{\mathscr{Z} , \geq h}{\longleftrightarrow}  \bar{\Gamma}     \big\} \cap \big\{   \mathscr{X} \overset{\mathscr{Z} , \geq h}{\not\longleftrightarrow}      \mathscr{Y}          \big\} \cap \big\{   \mathscr{X} \overset{\mathscr{Z} , \geq h}{\not\longleftrightarrow}   \Gamma   \big\} \cap \big\{ \mathscr{Y} \overset{\mathscr{Z} , \geq h}{\not\longleftrightarrow}        \Gamma \big\}    \big]          \\
       \equiv  \textbf{P}_h \big[ \mathcal{E}_1 \big] \textbf{P}_h \big[ \mathscr{E}_1 \cap \mathscr{E}_2 \cap \mathscr{E}_3 \cap \mathscr{E}_4 \cap \mathscr{E}_5 \big]  \\  \overset{(\mathrm{FKG})}{\leq} (\bar{\Gamma} \textit{Probability II}) \equiv  \textbf{P}_h \big[ \mathcal{E}_1 \cap \big\{ \mathscr{E}_1  \cap \mathscr{E}_2 \cap \mathscr{E}_3     \cap \mathscr{E}_4 \cap \mathscr{E}_5 \big\}                        \big]  \\ \leq    \textbf{P}_h \big[ \mathscr{X} \overset{\mathscr{Z},\geq h}{\longleftrightarrow}    \bar{\Gamma}   \big]     \textbf{P}_h \big[   \mathscr{Y} \overset{\mathscr{Z} , \geq h}{\longleftrightarrow} \bar{\Gamma}   \big]       \textbf{P}_h \big[    \mathscr{X} \overset{\mathscr{Z} , \geq h}{\not\longleftrightarrow}  \mathscr{Y} \big]         \textbf{P}_h \big[    \mathcal{E}_1    \big]                 \\   \overset{(\mathrm{FKG})}{\leq}   \textbf{P}_h \big[    \big\{   \big\{  \mathscr{X}           \overset{\mathscr{Z},\geq h}{\longleftrightarrow}   \bar{\Gamma} \big\} \cap \big\{   \mathscr{Y}             \overset{\mathscr{Z},\geq h}{\longleftrightarrow}           \bar{\Gamma}       \big\} \cap  \big\{           \mathscr{X}   \overset{ \mathscr{Z} , \geq h}{\not\longleftrightarrow}     \mathscr{Y}        \big\} \big\}   \cap \mathcal{E}_1       \big]             \\ \equiv  (\bar{\Gamma} \textit{Probability I})         \text{ } \text{ , }
       \end{align*}
       
       \noindent because of the fact that,
       
            \begin{align*}
        0< \textbf{P}_h \big[      \mathscr{E}_4                      \big]  < 1 \text{, }   \\    0< \textbf{P}_h \big[   \mathscr{E}_5                    \big]  < 1         \text{ } \text{ . }
       \end{align*}

     \noindent Therefore, it suffices to construct a map,

     \begin{align*}
        f_h :       \big\{    \mathcal{E}_1 \cap \mathscr{E}_1  \cap \mathscr{E}_2 \cap \mathscr{E}_3     \cap \mathscr{E}_4 \cap \mathscr{E}_5      \big\}  \longrightarrow      \big\{ \mathscr{X} \overset{\mathscr{Z} , \geq h }{\longleftrightarrow} \mathscr{Y}        \big\}                     \text{ } \text{ , }
     \end{align*}
     
     \noindent in which the pullback under $f^{-1}_h$ denotes the events contained within the intersection $ \big\{    \mathcal{E}_1 \cap \mathscr{E}_1  \cap \mathscr{E}_2 \cap \mathscr{E}_3     \cap \mathscr{E}_4 \cap \mathscr{E}_5      \big\}$. For some $\omega$ in the preimage space, and the pushforward of $\omega$, $f_h(\omega)$, under $f_h ( \cdot )$, to demonstrate that such an inequality holds, $f_h$ should satisfy,

     \begin{itemize}
         \item[$\bullet$]      \underline{Property 1}: The intersection between $\omega$, and the pushforward of $\omega$ under $f_h$, satisfies the following constraint,

    \begin{align*}
         \big| \omega \cap f_h (\omega) \big| < \mathscr{D}                           \text{, } 
    \end{align*}      
    
     for some finite $\mathscr{D}$, in which the intersection between the configuration $\omega$ and its pushforward, $f_h(\omega)$, can differ by at most $\mathscr{D}$ edges. 
         \item[$\bullet$]    \underline{Property 2}: Building off \underline{Property 1}, a configuration $\omega$ from the preimage space, and the pushforward under $f_h( \cdot )$, $f_h(\omega)$, satisfies the following constraint,

    \begin{align*}
     \big\{  \omega ,    \omega^{\prime} \in \{ \mathcal{E}_1 \cap \mathscr{E}_1 \cap \mathscr{E}_2 \cap \mathscr{E}_3 \cap \mathscr{E}_4 \cap \mathscr{E}_5 \big\}   :     \big|   f ( \omega ) \equiv f (\omega^{\prime})    \big|      < \mathscr{D}                     \big\}                                              \text{, } 
    \end{align*}    
         
       \noindent for two configurations $\omega$ and $\omega^{\prime}$ from the preimage space of $f_h$, in which it is stipulated that at most $\mathscr{D}$ edges between $f(\omega)$ and $f(\omega^{\prime})$ are equal, in which equality amongst edges of the pushforward configurations $f_h(\omega)$ and $f_h(\omega^{\prime})$ is given by,
       
       \begin{align*}
             \big|   f_h ( \omega ) \equiv f_h (\omega^{\prime})    \big|     = \mathscr{D}     \Longleftrightarrow  \forall   e_1 (\omega) \in f_h(\omega) , \cdots ,   e_n(\omega)  \in f_h( \omega) , e_1 (\omega^{\prime}) \in f_h(\omega^{\prime}), \\  \cdots , e_n ( \omega^{\prime} )\in f_h(\omega^{\prime}) , \exists \mathscr{D} < n  :  e_1 ( \omega ) \equiv e_1 ( \omega^{\prime} ) , \cdots , e_{\mathscr{D}} ( \omega ) \equiv e_{\mathscr{D}} ( \omega^{\prime} )   \text{, } 
       \end{align*}
       
       \noindent in which there is equality amongst the first $\mathscr{D} < n$ edges associated with the pushfoward configurations $f_h(\omega)$ and $f_h(\omega^{\prime})$.
       
     \end{itemize}
     
     \noindent With regards to \underline{Property 1} and \underline{Property 2}, of each condition holds, observe that ($\bar{\Gamma}$ \textit{Probability I}) would hold with the constant in the upper bound taking the form,
     
     \begin{align*}
         \mathcal{C}_1 \equiv  \mathscr{D} \big(     \mathrm{inf} \big\{        0  , 1 - \delta   \big\}    \big)^{-\mathscr{D}}           \text{. }
     \end{align*}
     
     \noindent for suitable $\delta>0$. To demonstrate that \underline{Property 1} holds, first, introduce, for some path of strictly positive length,
     
       \begin{align*}
          \pi_{u,n} \in \pi_{\mathscr{U},\mathscr{N}} \equiv  \underset{n \in \mathscr{N} \subset \textbf{N}}{\underset{u \in \mathscr{U}}{\bigcup}} \pi_{u,n} \equiv \underset{n \in \mathscr{N} \subset \textbf{N}}{\underset{u \in \mathscr{U}}{\bigcup}} \bigg\{     \forall u \in \mathscr{U} , n \in \mathscr{N}, \exists          u_1 , \cdots , u_{\mathscr{D}} , \cdots , u_n \in \mathcal{V}           :    \big\{ u_1 \equiv u_n  \big\} \cup  \cdots \\     \big\{  \underset{0\leq i\leq n}{\bigcup} \big\{  \pi_{u_i}              \overset{ \mathscr{Z} , < h}{\longleftrightarrow}   \pi_{u_{i+1}}   \big\}    \big\} \bigg\}  \text{, } 
     \end{align*}
     
     \noindent where $u$ are drawn from the set, for $\mathscr{X}^{\prime}\supsetneq \mathscr{X}$,

     \begin{align*}
   \mathscr{U} \equiv \big\{    u \in \bar{\Gamma}      :     \textbf{P}_h \big[    u       \overset{\mathscr{Z},< h }{\longleftrightarrow}               \mathscr{X} \big] > 0 , \textbf{P}_h \big[     u        \overset{\mathscr{Z}, \mathscr{X},< h }{\longleftrightarrow}     \mathscr{X}^{\prime}          \big] > 0  ,   \textbf{P}_h \big[    \big\{   u                           \overset{\mathscr{Z},< h }{\longleftrightarrow}    \mathscr{X}  \big\} \cap \big\{    u                 \overset{\mathscr{Z},< h }{\longleftrightarrow}    \mathscr{X}^{\prime}    \big\} \cap \big\{  \mathscr{X}^{\prime}    \overset{\mathscr{Z},< h }{\longleftrightarrow}        u       \big\}             \big] = 0  \big\}       \text{. }
     \end{align*}
     
     \noindent Similarly, for another point $v$ instead of $u$, the self-avoiding path from the level-set connectivity event of the metric-graph GFF would instead be dependent upon the conditions,
     
     \begin{align*}
\pi_{v , n^{\prime}}  \propto  \big\{   v_1 \equiv v_{n^{\prime}}     \big\}      \cup \bigg\{  \underset{0\leq i\leq n^{\prime}}{\bigcup} \big\{  \pi_{v_i}              \overset{ \mathscr{Z} , < h}{\longleftrightarrow}   \pi_{v_{i+1}}   \big\}    \bigg\}          \text{, } 
     \end{align*}
     
     \noindent being fulfilled. 
     
     \bigskip
     
  \noindent \underline{Case A}

     \bigskip
     
     \noindent From $\pi_{u,n}$ and $\pi_{v,n^{\prime}}$, as two self-avoiding paths associated with the existence of a level-set percolation crossing in the metric-graph GFF, to demonstrate that \underline{Property 1} and \underline{Property 2} of $f_h( \cdot)$ hold, fix some $u \in \mathscr{U}$, $v \in \mathscr{V}$, where, as defined with $\mathscr{U}$ above, introduce, for $\mathscr{V}^{\prime}\supsetneq \mathscr{V}$,

     \begin{align*}
          \mathscr{V} \equiv \big\{    v \in \bar{\Gamma}      :     \textbf{P}_h \big[    v       \overset{\mathscr{Z},< h }{\longleftrightarrow}               \mathscr{V} \big] > 0 , \textbf{P}_h \big[ v        \overset{\mathscr{Z}, \mathscr{V},< h }{\longleftrightarrow}     \mathscr{V}^{\prime}          \big] > 0  ,   \textbf{P}_h \big[    \big\{   v                          \overset{\mathscr{Z},< h }{\longleftrightarrow}    \mathscr{V}  \big\} \cap \big\{    v                 \overset{\mathscr{Z},< h }{\longleftrightarrow}    \mathscr{X}^{\prime}    \big\} \cap \big\{  \mathscr{V}^{\prime}    \overset{\mathscr{Z},< h }{\longleftrightarrow}        v       \big\}             \big] = 0  \big\}       \text{, } 
     \end{align*}

     \noindent corresponding to $v \in \bar{\Gamma}$ satisfying the three conditions provided above. Pertaining to the pushforward of $\omega$ in the preimage space under $f_h(\cdot)$ fix $u,v$ so that $\overline{\{u \}} \equiv \overline{\{ v \}}$. Under this assumption, $f_h(\omega)$ is modified from $\omega$ in the following ways, which we denote below with \underline{Case A}.
     
     \bigskip

\begin{itemize}
    \item[$\bullet$] \underline{Modification one} (\textit{metric-graph GFF analog for closing edges in configurations as previously argued for Bernoulli percolation}): For edges in $f(\omega)$ for which the event,
    
    \begin{align*}
    \textbf{P}_h \big[          \mathscr{X} \overset{ \mathscr{Z} ,\geq h}{\longleftrightarrow}  \mathscr{Y}         \big]  \text{, } 
    \end{align*}

    \noindent occurs with positive probability, close all edges belonging to the following set, which satisfy the requirements,

    \begin{align*}
                \mathscr{V}^{\prime} \equiv \big\{ \forall v  \in \mathcal{V} ,  \exists v \in f_h( \omega) : v \cap \overline{\{ u \}} \neq \emptyset \big\}                   \text{, } 
    \end{align*}
    
    \noindent corresponding to the set of edges belonging to $f_h(\omega)$ which have nonempty intersection with $\overline{\{ u \}}$, as well as the set,
  
    \begin{align*}
            \mathscr{V}^{\prime}_{\pi_{u}} \equiv \mathscr{V}^{\prime}_u        \equiv \big\{      \forall u \in \mathcal{V} , \exists  u \in \pi_u :         u \cap \overline{\{ u \}} \neq \emptyset      \big\}     \text{, } 
    \end{align*}

    \noindent corresponding to the set of edges belonging to $\pi_u$ for which the intersection of $\pi_u$ with $\overline{\{ u\}}$ is nonempty. From $\mathscr{V}^{\prime}$ and $\mathscr{V}^{\prime}_u$, admissible $u$ and $v$ satisfying conditions from the intersection,

    \begin{align*}
 \mathscr{V}^{\prime}  \cap  \mathscr{V}^{\prime}_u     \text{, } 
    \end{align*}
    
    \noindent correspond to vertices from the vertex set $\mathcal{V}$ of the \textit{three-dimensional environment} that, besides each having nonempty intersection with $\overline{\{ u\}}$, also lie within the intersection of vertices of $f_h(\omega)$ and $\pi_u$, for which,
    
    \begin{align*}
     \textbf{P}_h \big[          v     \overset{\mathscr{Z} , < h}{\longleftrightarrow}      \overline{\{ u \}}        \big] > 0    \text{, } 
    \end{align*}
    
    \noindent for $v \in \mathscr{V}^{\prime}$, and also for which,

    \begin{align*}
        \textbf{P}_h \big[   u       \overset{\mathscr{Z} , < h}{\longleftrightarrow}                    \overline{\{ u \}}       \big] > 0       \text{, } 
    \end{align*}
    
    \noindent for $u \in \mathscr{V}^{\prime}_u$. From the modifications introduced for \underline{Modification one}, the two probabilities above would instead have vertices, and hence edges, of the pushforward configuration $f_h(\omega)$ satisfy,
    
    \begin{align*}
        \textbf{P}_h \big[          v     \overset{\mathscr{Z} , \geq  h}{\longleftrightarrow}      \overline{\{ u \}}        \big] > 0    \text{, } 
    \end{align*}

    \noindent and also,
    
     \begin{align*}
           \textbf{P}_h \big[   u       \overset{\mathscr{Z} , \geq  h}{\longleftrightarrow}                    \overline{\{ u \}}       \big] > 0  \text{, } 
    \end{align*}
    
    \noindent for edges satisfying,

    \begin{align*}
   \underset{ e(v_i,v_{i+1}) \overset{i \equiv  n-1}{=}  e ( v_{n-1} , \overline{\{ u\}} )   }{\underset{e(v,v_i) \overset{i \equiv 1}{=} e ( v, v_1) }{\underset{\forall  1 \leq  i     \leq n          }{\underbrace{\textbf{1}_{e(v , v_i) \overset{\mathscr{Z} , \geq h}{\longleftrightarrow}    e(v_i , v_{i+1})           }              \equiv 1 } }}} \propto \underset{v_j \overset{j\equiv \frac{n}{2}}{=}    \text{ }   \overline{\{ u \}}      }{\underset{v_j \overset{j \equiv 1}{=}  v    }{\underset{\forall  1 \leq j \leq \frac{n}{2}                }{\underbrace{\textbf{1}_{ v_j       \overset{\mathscr{Z} , \geq h}{\longleftrightarrow}    v_{j+1}      }     \equiv 1    }}}}   \text{, }   \text{ }              \text{. }
    \end{align*}
    
    \noindent corresponding to the event $\{    v \overset{\mathscr{Z} , \geq h}{\longleftrightarrow} \overline{\{ u\}}   \}$, while similarly, the indicator for each index is satisfied for,
    
    \begin{align*}
    \textbf{1}_{e ( u , u_j )\overset{\mathscr{Z} , \geq h} {\longleftrightarrow}        e(u_j , u_{j+1} )        } \equiv 1                 \text{, } 
    \end{align*}

    \noindent $\forall 1 \leq j \leq n^{\prime} $,corresponding to the event $\{     u \overset{\mathscr{Z} , \geq h}{\longleftrightarrow} \overline{\{ u\}}       \}$, under a similar convention of indices as provided for $\{    v \overset{\mathscr{Z} , \geq h}{\longleftrightarrow} \overline{\{ u\}}   \}$.

    \item[$\bullet$] \underline{Modification two} (\textit{metric-graph GFF analog for opening line segments in three dimensions as previously argued for Bernoulli percolation}). In opposition to \underline{Modification 1} for closing edges for metric-graph GFF level-set percolation events, for \underline{Modification two}, we consider the following procedure for opening edges depending upon the height of the metric-graph GFF, in which, from the edge set,
    
    \begin{align*}
            \mathcal{E}    \supsetneq \mathscr{E}^{\prime\prime} \equiv \big\{    \forall e \in \mathcal{E} , \exists e \in \mathcal{E} \big( \overline{\{ u \}} \big) :     e \cap \gamma \neq \emptyset                            \big\}           \text{, } 
    \end{align*}
    
    \noindent where,
    
    \begin{align*}
           \big| \gamma \big(    u   ,     \overline{\big\{ u \big\}}     \big) \big|    \equiv   \big| \gamma \big|           \equiv  \underset{\Gamma \cap u \cap \overline{\{ u \}} \neq \emptyset }{\underset{{\Gamma} \subset \mathcal{S}}{\mathrm{inf}} }      \big\{      \rho_{\Gamma} \big(         u                                 ,            \overline{\big\{ u \big\}}     \big)      \big\}      \text{, } 
    \end{align*}

    \noindent with edge set from $\mathscr{E}^{\prime\prime}$ given by,

    \begin{align*}
  \mathcal{E} \big( \overline{\{ u \}} \big) \equiv \big\{   e \in \mathcal{E} : e \cap \overline{\{ u\}} \neq \emptyset    \big\}       \text{, } 
    \end{align*}

    \noindent for the restriction of the metric $\rho$ to $\Gamma \subset \mathcal{S}$, which is defined by,

    \begin{align*}
          \rho_{\Gamma} \big( u , \overline{\big\{ u \big\} }  \big) \equiv   \big\{ \forall v_1 , v_2 \in \mathcal{V} , \exists v_1 \equiv u , v_2 \equiv \overline{\big\{ u \big\}} :  \rho_{\Gamma}   \big( u , \overline{\big\{ u \big\} }   \big)     \equiv \rho_{\mathcal{S}|_{\Gamma}} \big(   u , \overline{\big\{ u \big\} }   \big)     \big\}        \text{ }\text{ , } 
    \end{align*}
    
    \noindent and,

    \begin{align*}
   \mathcal{E}   \supsetneq \mathcal{E} \big( \overline{\{ u \}} \big) \equiv \big\{    e \in \mathcal{E} : e \cap \overline{\big\{ u \big\}    } \neq \emptyset \big\}  \text{, } 
    \end{align*}

    \noindent as all edges belonging to $\mathscr{E}^{\prime\prime}$ that belong to the path indicated above are opened.

    \item[$\bullet$] \underline{Modification three} (\textit{opening edges from the shortest path emanating from $v$, in comparison to opening edges from the shortest path emanating from $u$ as provided in} \underline{Modification two}). For the final modification, to build off arguments introduced for \underline{Modification two}, introduce the following edge set,

    \begin{align*}
   \mathscr{E}^{\prime\prime\prime}  \big( v , e_2 \big) \equiv \mathscr{E}^{\prime\prime\prime} \equiv \big\{   \forall e_1 \in \mathcal{E} \big( v \big)  , \exists                     e_2 \in \mathcal{E} \big( \overline{\big\{ u \big\}} \big) :                          \Gamma^{\prime}     \big( e_1 , e_2 \big)        \cap \rho_{\Gamma} \neq \emptyset      \big\}  \text{, } 
    \end{align*}

    \noindent where,

    \begin{align*}
    \mathcal{E} \big( v \big) \equiv     \big\{ e \in \mathcal{E} : e \cap v \neq \emptyset   \big\}          \text{, } 
    \end{align*}

    \noindent and,

    \begin{align*}
           \Gamma \big(    e_1     ,     e_2      \big) \equiv   \big\{  \overline{\big\{ u \big\}} \in \mathcal{E}    : \Gamma^{\prime} \big( e_1 , e_2 \big) \cap \overline{\big\{ u \big\}} \neq \emptyset     \big\}              \text{, } 
    \end{align*}

    \noindent with the length of the path given by the infimum,

    \begin{align*}
   \big| \Gamma^{\prime} ( e_1 , e_2 ) \big| \equiv \big| \Gamma^{\prime} \big| \equiv    \underset{\Gamma \subset \mathcal{S}}{\underset{e_1 , e_2 \in \mathcal{E}}{ \mathrm{inf}}}      \big\{             \rho_{\Gamma} \big(   e_1   ,  e_2  \big)      \big\}       \text{, } 
    \end{align*}

    \noindent in which we incorporate conditions from $\gamma$ introduced in \underline{Modification two}, namely that this shortest path, with respect to the infimum of the metric $\rho$, intersect $v$ and $\overline{\big\{ u \big\}}$. Concluding, opening the edges in $\overline{\big\{ u \big\}}$ with this property entails that,
    
    \begin{align*}
         \textbf{P}_h \big[           e_1       \overset{\mathscr{Z} ,       < h  }{\longleftrightarrow}  \big(      \Gamma^{\prime} \big( e_1 , e_2 \big)   \cup       \gamma     \big)       \big]       \text{, } 
    \end{align*}
    
    \noindent occurs with positive probability, in addition to,

    \begin{align*}
          \underset{\Gamma^{\prime} \subset \mathcal{S}}{\underset{\gamma \subset \mathcal{S}}{\mathrm{inf}}} \big\{    \textbf{P}_h \big[           e_1       \overset{\mathscr{Z} ,       < h  }{\longleftrightarrow}  \big(      \Gamma^{\prime} \big( e_1 , e_2 \big)   \cup       \gamma     \big)       \big]    \big\}       \text{, } 
    \end{align*}

    \noindent occurring with positive probability, for an infimum taken over $\gamma, \Gamma^{\prime} \subset \mathcal{S}$.

\end{itemize}

       \noindent From \underline{Modification one}, \underline{Modification two}, and \underline{Modification three} described above, observe,
      
        \begin{align*}
                \big| \omega \cap  f_h \big( \omega \big) \big| < \mathscr{D} \text{, } 
     \end{align*}
     
     \noindent as the intersection between the number of common edges between $\omega$ and $f_h\big( \omega \big)$ can be upper bounded by a suitable $\mathscr{D}$ as provided before describing the three modifications of $f_h(\omega)$. As a result, upon performing the respective modifications of closing, opening, and opening, certain edges depending upon whether the level-set percolation event across height $\geq h$ or $< h$ occurs, $\overline{\big\{ u \big\}}$ can be uniquely reconstructed in $f_h(\omega)$, in which,
     
    \begin{align*}
        \textbf{P}_h \big[  \big\{  \overline{\big\{ z \big\} }       \overset{\mathscr{Z} , \geq h}{\longleftrightarrow}         \mathscr{X} \big\} \cup   \big\{      \overline{\big\{ z \big\} }     \overset{\mathscr{Z} , \geq h}{\longleftrightarrow}    \mathscr{Y}   \big\}        \big]   > 0     \text{, } 
    \end{align*}
     
       \noindent as the two events displayed above, namely the union of $\overline{\big\{ z \big\}}$ being connected to $\mathscr{X}$, in addition to $\overline{\big\{ z \big\}}$ being connected to $\mathscr{Y}$, occurs with positive probability. 
       
       \bigskip
       
       \noindent  \underline{Case B}

       \bigskip
       
       \noindent On the other hand, if $\overline{\mathscr{U}} \cap \overline{\mathscr{V}} \equiv \emptyset$, then a similar series of modifications, as given previously, can be applied to the pushfoward $f_h(\omega)$. For the following items, fix $u \in \mathscr{U}$, and $v \in \mathscr{V}$. From such a choice of $u$ and $v$, one has,
       
       \begin{align*}
          \textbf{P}_h \big[            \overline{\big\{ u \big\}} \overset{\mathscr{Z} , \geq h}{\longleftrightarrow}        \mathscr{Y}      \big] \text{, } 
       \end{align*}
       
       \noindent occurring with zero probability, in addition to,
       
           \begin{align*}
           \textbf{P}_h \big[           \overline{\big\{ v\big\}}      \overset{\mathscr{Z} , \geq h}{\longleftrightarrow} \mathscr{X}   \big]   \text{, } 
       \end{align*}
       
         \noindent occurring with zero probability.

     \begin{itemize}
         \item[$\bullet$] \underline{Modification one} (\textit{closing edges, as an alternative modification to} $f_h\big(\omega\big)$ \textit{presented previously in Modification one}). Introduce the set of edges,
         
         \begin{align*}
         \mathscr{E}^{\prime\prime\prime\prime} \equiv \big\{   e \in \mathcal{E} : e \cap \big( \overline{\big\{ u \big\}} \cup \overline{\big\{ v \big\}} \big) \neq \emptyset                           \big\}    \text{, } 
         \end{align*}
         
      \noindent in which one examines the edges, of the \textit{three-dimensional environment} that have nonempty intersection with $\overline{\big\{ u\big\}} \cup \overline{\big\{ v \big\}}$. Furthermore, related objects for this modification to $f_h \big( \omega \big)$ also include,
      
      \begin{align*}
        \mathscr{E}^{\prime\prime\prime\prime\prime}_{u, < h} \equiv \mathscr{E}_{u , < h} \equiv \bigg\{         e \in \mathcal{E} : \big(  e \cap \pi_u \neq \emptyset    \big) \cup \bigg(    \textbf{P}_h \big[       \big(   \overline{\big\{ u \big\}} \cup \overline{\big\{ v \big\}}  \big)  \overset{\mathscr{Z}, < h }{\longleftrightarrow}   \pi_u               \big]                \bigg)     \bigg\} \text{, } 
      \end{align*}
      
      \noindent corresponding to the edges of the \textit{three-dimensional environment} which have nonempty intersection with $\pi_u$, and,
      
      \begin{align*}
         \mathscr{E}^{\prime\prime\prime\prime\prime\prime}_{v,<h} \equiv       \mathscr{E}_{v,<h}   \equiv  \bigg\{    e \in \mathcal{E} : \big(  e \cap \pi_v \neq \emptyset   \big) \cup \bigg(       \textbf{P}_h \big[    \big(   \overline{\big\{ u \big\}} \cup \overline{\big\{ v \big\}}  \big)  \overset{\mathscr{Z}, < h }{\longleftrightarrow}           \pi_v      \big]           \bigg)                \bigg\}         \text{, } 
      \end{align*}
      
      \noindent corresponding to the edges of the \textit{three-dimensional environment} which have nonempty intersection with $\pi_v$, and, 
      
      \begin{align*}
         \mathscr{E}^{\prime\prime\prime\prime\prime\prime\prime}_{\Gamma,<h} \equiv  \mathscr{E}_{\Gamma, <h}    \equiv  \bigg\{ e  \in \mathcal{E} : \big(  e \cap \Gamma \neq \emptyset   \big) \cup \bigg(          \textbf{P}_h \big[            \big(  \overline{\big\{ u \big\}} \cup \overline{\big\{ v \big\}}  \big)  \overset{\mathscr{Z}, < h }{\longleftrightarrow}   \Gamma     \big]                 \bigg)               \bigg\}      \text{, } 
      \end{align*}

      \noindent corresponding to the edges of the \textit{three-dimensional environment} which From $ \mathscr{E}^{\prime\prime\prime\prime}$, $\mathscr{E}^{\prime\prime\prime\prime\prime}_u$, $ \mathscr{E}^{\prime\prime\prime\prime\prime\prime}_v$, and $  \mathscr{E}^{\prime\prime\prime\prime\prime\prime\prime}_{\Gamma}$ introduced above, observe that the series of modifications,
      
      \begin{align*}
          \textbf{P}_h \big[   e \in \mathcal{E} \cap  \mathscr{E}^{\prime\prime\prime\prime} \cap \mathscr{E}^{\prime\prime\prime\prime\prime}_{u, < h} \cap  \mathscr{E}^{\prime\prime\prime\prime\prime\prime}_{v, < h} \cap            \mathscr{E}^{\prime\prime\prime\prime\prime\prime\prime}_{\Gamma, < h}             \big]   > 0       \overset{\mathrm{'edge\text{ }  closing'}}{\Longleftrightarrow}          \textbf{P}_h \big[    e \in \mathcal{E} \cap  \mathscr{E}^{\prime\prime\prime\prime} \cap \mathscr{E}^{\prime\prime\prime\prime\prime}_{u, \geq h} \cap  \mathscr{E}^{\prime\prime\prime\prime\prime\prime}_{v, \geq h} \cap            \mathscr{E}^{\prime\prime\prime\prime\prime\prime\prime}_{\Gamma, \geq h}             \big]  > 0        \text{, } 
      \end{align*}
      
      \noindent incorporate changes from open edges to closed edges, depending upon whether the metric-graph GFF attains a sufficiently high level relative to the height threshold $h$. 
      
      \item[$\bullet$] \underline{Modification two} (\textit{opening edges from the shortest path between} $\Gamma$ \textit{and} $u$). Introduce the set of edges,
      
     \begin{align*}
             \mathscr{E}^{\prime\prime\prime\prime\prime\prime\prime\prime}_{u , \Gamma, \geq h} \equiv \mathscr{E}_{u , \Gamma, \geq h} \equiv   \big\{ e \in \mathcal{E} :                  \big\{  e \cap \big\{ u \overset{\mathscr{Z} , \geq h}{\longleftrightarrow} \Gamma \big\}  \neq \emptyset \big\} \cup \big\{ d_{u , \Gamma}  \big\}             \big\}          \text{, } 
     \end{align*} 
      
      \noindent where $d_{u,\Gamma}$ satisfies,
      
      \begin{align*}
         d_{u , \Gamma}   \equiv            {\underset{ u \in \mathscr{U}}{\underset{\Gamma , \Gamma^{\prime} \subset \mathcal{S}}{\mathrm{inf}} }}    \big\{     \Gamma^{\prime} \big(    \Gamma  , u  , \big\{ \Gamma  \overset{\mathscr{Z} , \geq h}{\longleftrightarrow}  u  \big\}       \big)                                             \big\}          \text{, } 
      \end{align*}

      \noindent with $\Gamma^{\prime}$ satisfying the three following conditions, simultaneously,
      
      \begin{align*}
        \Gamma^{\prime} \equiv \big\{     \forall     e \in \mathcal{E} , \exists n^{\prime\prime} \in \textbf{N} : \big\{  \big\{  e \cap  \Gamma \big\} \neq \emptyset \big\} \cup \big\{    \big\{   u \cap \Gamma \big\} \neq \emptyset                      \big\} \cup    \big\{     d_{u, \Gamma} \equiv n^{\prime\prime}          \big\}                         \big\}   \text{. }
      \end{align*}
      
         \noindent As a result, for suitable paths $\gamma^{\prime}$,
         
         \begin{align*}
         \textbf{P}_h \big[                 \big\{ e \in    \mathscr{E}^{\prime\prime\prime\prime\prime\prime\prime\prime}_{u , \Gamma , \geq h}  \big\} \cup \big\{ \gamma^{\prime} \equiv \Gamma^{\prime} \big( \Gamma , u , \big\{   \Gamma \overset{\mathscr{Z} , \geq h}{\longleftrightarrow}  u \big\} \big)  \big\}   \big]   \overset{\mathrm{'edge\text{ }  opening'}}{\Longleftrightarrow}  \textbf{P}_h \big[      \big\{ e \in    \mathscr{E}^{\prime\prime\prime\prime\prime\prime\prime\prime}_{u , \Gamma , < h}  \big\} \cup \cdots \\  \big\{ \gamma^{\prime} \equiv \Gamma^{\prime}   \big( \Gamma , u , \big\{   \Gamma \overset{\mathscr{Z} , < h}{\longleftrightarrow}  u \big\} \big)\big\}                  \big]  \text{ , } 
         \end{align*}
         
         \noindent in which, relative to the height threshold crossing $h$, edges in the pushforward configuration $f_h \big( \omega \big)$ are opened, along the shortest such path $\Gamma^{\prime}$.
         
         \item[$\bullet$] \underline{Modification three} (\textit{opening edges from the shortest path between} $\Gamma$ \textit{and} $v$). Lastly, introduce the set of edges,
         
         \begin{align*}
             \mathscr{E}^{\prime\prime\prime\prime\prime\prime\prime\prime\prime}_{v , \Gamma , \geq h } \equiv \mathscr{E}_{v , \Gamma , \geq h}   \equiv \big\{              e \in \mathcal{E} :                       \big\{  e \cap \big\{ v \overset{\mathscr{Z} , \geq h}{\longleftrightarrow} \Gamma \big\}  \neq \emptyset \big\} \cup \big\{ d_{v , \Gamma}  \big\}                                          \big\}               \text{, } 
         \end{align*}

         \noindent  where $d_{v,\Gamma}$ satisfies, as does $d_{u,\Gamma}$,
      
      \begin{align*}
         d_{u , \Gamma}   \equiv            {\underset{ v \in \mathscr{V}}{\underset{\Gamma , \Gamma^{\prime\prime} \subset \mathcal{S}}{\mathrm{inf}} }}    \big\{     \Gamma^{\prime} \big(    \Gamma  , v  , \big\{ \Gamma  \overset{\mathscr{Z} , \geq h}{\longleftrightarrow}  v  \big\}       \big)                                             \big\}          \text{, } 
      \end{align*}
         
        \noindent with $\Gamma^{\prime}$ satisfying the three following conditions, simultaneously,
      
      \begin{align*}
        \Gamma^{\prime\prime} \equiv \big\{     \forall     e \in \mathcal{E} , \exists n^{\prime\prime\prime} \in \textbf{N} : \big\{  \big\{  e \cap  \Gamma \big\} \neq \emptyset \big\} \cup \big\{    \big\{   v \cap \Gamma \big\} \neq \emptyset                      \big\} \cup    \big\{     d_{v, \Gamma} \equiv n^{\prime\prime\prime}          \big\}                         \big\}   \text{. }
      \end{align*}
      
         \noindent  As a result, for suitable paths $\gamma^{\prime}$,
         
         \begin{align*}
         \textbf{P}_h \big[                 \big\{ e \in    \mathscr{E}^{\prime\prime\prime\prime\prime\prime\prime\prime\prime}_{v , \Gamma , \geq h}  \big\} \cup \big\{ \gamma^{\prime} \equiv \Gamma^{\prime\prime} \big\}   \big]   \overset{\mathrm{'edge\text{ }  opening'}}{\Longleftrightarrow}  \textbf{P}_h \big[      \big\{ e \in    \mathscr{E}^{\prime\prime\prime\prime\prime\prime\prime\prime\prime}_{v , \Gamma , < h}  \big\} \cup \big\{ \gamma^{\prime} \equiv \Gamma^{\prime\prime} \big\}                  \big] \text{, } 
         \end{align*}

        \noindent in which, as previously argued for \underline{Modification two}, implies, relative to the height threshold crossing $h$, that edges in the pushforward $f_h(\omega)$ can have edges opened as demonstrated in the fashion above.
         
     \end{itemize}

       \noindent As with \underline{Modification one}, \underline{Modification two}, and \underline{Modification three} introduced for the case $\overline{\big\{ u \big\}} \equiv \overline{\big\{ v \big\}}$, if $\overline{\mathscr{U}} \cap \overline{\mathscr{V}} \equiv \emptyset$, another strictly positive upper bound between the intersection of common edges between $\omega$ and $f_h\big( \omega \big)$ satisfies,
       
       \begin{align*}
            \big|  \omega \cap f_h \big( \omega \big) \big| < \mathscr{D}^{\prime}     \text{ }\text{ . } 
       \end{align*}

       \noindent Hence, for $\overline{\mathscr{U}} \cap \overline{\mathscr{V}} \equiv \emptyset$, $\overline{\big\{ u \big\} }$ can be uniquely reconstructed from $f_h \big( \omega \big)$, in which,
       
       \begin{align*}
          \textbf{P}_h \big[ \mathscr{X}  \overset{ \Gamma^{\prime\prime} ,\geq h}{\longleftrightarrow}      \mathscr{Y} \big]                                    \text{. }
       \end{align*}
       
       \noindent Across the two subcases within \underline{Case A}, it has been demonstrated that each possibility holds with $\mathscr{D}$, and with $\mathscr{D}^{\prime}$, respectively. To proceed, consider the following, \underline{Case B}.

        \bigskip
 
       \noindent For the following analysis of the second case, we make use of a similar lower bound from the following intersection of events, for suitable, strictly positive, $\mathcal{C}_1$,

       \begin{align*}
              \textbf{P}_h \big[ \mathcal{E}_1 \cap \mathscr{E}_1  \cap \mathscr{E}_2 \cap \mathscr{E}_3     \cap \mathscr{E}_4 \cap \mathscr{E}^{\prime}_5                        \big] \leq                           \mathcal{C}_1            \textbf{P}_h \big[ \mathscr{X} \overset{\mathscr{Z} , \geq h}{\longleftrightarrow}   \mathscr{Y}          \big]                                   \text{, } \tag{$\bar{\Gamma}$ \textit{Probability III}}
       \end{align*}

       \noindent with the exception that the last event is,
    
        \begin{align*}
     \mathscr{E}^{\prime}_5  \equiv  \big\{   \mathscr{Y} \overset{\mathscr{Z},\geq h}{ \longleftrightarrow}  \Gamma  \big\}          \text{ } \text{ , }
        \end{align*}

       \noindent instead of,
       
        \begin{align*}
     \mathscr{E}_5  \equiv  \big\{   \mathscr{Y} \overset{\mathscr{Z},\geq h}{\not\longleftrightarrow}  \Gamma  \big\}         \text{ } \text{ . }
        \end{align*}
    
       \noindent In the arguments below, we very briefly summarize the points distinguishing \underline{Case C}, and \underline{Case D}, each of which respectively share in the arguments for \underline{Case A}, and for \underline{Case B}, in which edges are closed and/or opened.
       
       \bigskip

       \noindent As for \underline{Case A} and for \underline{Case B}, denote the map,

       \begin{align*}
                 f_h :       \big\{    \mathcal{E}_1 \cap \mathscr{E}_1  \cap \mathscr{E}_2 \cap \mathscr{E}_3     \cap \mathscr{E}_4 \cap \mathscr{E}^{\prime}_5      \big\}  \longrightarrow      \big\{ \mathscr{X} \overset{\mathscr{Z} , \geq h }{\longleftrightarrow} \mathscr{Y}        \big\}      \text{, } 
       \end{align*}
       
       \noindent in which we would like to consider the action of the pushforward $f_h \big( \cdot \big)$ on $\omega$, where,
       
       \begin{align*}
       \omega \in \big\{   \mathcal{E}_1 \cap \mathscr{E}_1  \cap \mathscr{E}_2 \cap \mathscr{E}_3     \cap \mathscr{E}_4 \cap \mathscr{E}^{\prime}_5      \big\}     \text{, } 
       \end{align*}
       
       \noindent is some configuration from the preimage space of $f_h \big( \cdot \big)$. For the two cases below, fix the same objects, as in previous cases, which include $u \in \mathscr{U}$, with $\big\{      \mathscr{Y}             \overset{\mathscr{Z} , \geq h }{\longleftrightarrow}          \Gamma   \big\} \supsetneq \big\{    \mathscr{Y}             \overset{\mathscr{Z} \backslash  \overline{\{ u \}} , \geq h }{\longleftrightarrow}          \Gamma       \big\}$.
       
       \bigskip
       
       \noindent \underline{Case C}
       
       \bigskip 
       
       \begin{itemize}
       \item[$\bullet$] \underline{Modification one} (\textit{closing edges of the metric-graph GFF}).  
       \noindent The sequence of arguments slightly differs from that presented in \underline{Case A}. Namely, from the event given in the preimage space, introduce the set of edges,
       
       \begin{align*}
        \mathscr{E}^6 \equiv         \mathscr{E}_{\overline{\{ u\}}, \geq h}   \equiv \big\{    \forall   e \in \mathcal{E} , \exists    e \cap \overline{\big\{ u \big\}} \neq \emptyset   :                \textbf{P}_h \big[      \big\{   \Gamma                          \overset{\mathscr{Z} , \geq h}{\longleftrightarrow}       \overline{\{ u \}}        \big\}      \cup            \big\{   \Gamma \cap e \equiv \emptyset   \big\}  \big]  > 0        \big\}    \text{, } 
       \end{align*}

       \noindent corresponding to the edges of the \textit{three-dimensional environment} with nonempty intersection with $\overline{\big\{ u \big\}}$. Equipped with $\mathscr{E}_6$, additionally introduce,

       \begin{align*}
         \big(  \mathscr{E}^6 \big)^{\prime} \equiv \mathscr{E}^{\prime}_{  \overline{\{ u\}}, \geq h}  \equiv \big\{   \forall     e \in \mathcal{E} , \exists          e \cap         \overline{\{ u \}}                 \neq \emptyset   :        e \cap \overline{\{ u \}} \cap  \pi_u \equiv \pi_u    , \textbf{P}_h \big[              \overline{\{ u\}}        \overset{\mathscr{Z} , \geq h}{\longleftrightarrow}         \pi_u               \big]  \equiv 0    \big\}            \text{, } 
       \end{align*}

       \noindent and,
       
          \begin{align*}
               \big(\mathscr{E}^6 \big)^{\prime\prime} \equiv  \mathscr{E}^{\prime\prime}_{ \overline{\{ u\}}, \geq h} \equiv  \big\{ \forall     e \in \mathcal{E} , \exists          e \cap        \overline{\{ u \}}                  \neq \emptyset   : e \cap  \overline{\{ u \}}   \cap \Gamma \equiv \Gamma  ,  \textbf{P}_h \big[      \overline{\{ u \} }        \overset{\mathscr{Z} , \geq h}{\longleftrightarrow}   \Gamma                   \big]  \equiv 0  \big\}            \text{. }
       \end{align*}
       
       \noindent Concluding, performing the 'edge closing' procedure in from the edge constraints of $\mathcal{E}$ provided from $\mathscr{E}_6$, $\mathscr{E}^{\prime}_6$ and $\mathscr{E}^{\prime\prime}_6$ yields,

       \begin{align*}
          \textbf{P}_h \big[             e \in \mathscr{E}_{ \overline{\{ u\}} , < h  } \cap \mathscr{E}^{\prime}_{ \overline{\{ u\}} , < h  }  \cap \mathscr{E}^{\prime\prime}_{  \overline{\{ u\}} , < h }                        \big] \\  \overset{\mathrm{'edge\text{ }  closing'}}{\Longleftrightarrow}    \\ \textbf{P}_h \big[     e \in \mathscr{E}_{ \overline{\{ u\}} , \geq h }  \cap \mathscr{E}^{\prime}_{ \overline{\{ u\}} , \geq h } \cap \mathscr{E}^{\prime\prime}_{  \overline{\{ u\}} ,  \geq h  }                             \big]  \text{ } \text{ . }
       \end{align*}

       \item[$\bullet$] \underline{Modification two} (\textit{opening edges of the metric-graph GFF after closing edges of the metric-graph GFF}). Introduce the set of edges,
       
       \begin{align*}
            \mathscr{E}_7 \equiv \mathscr{E}_7 \big(       e  , \Gamma      ,   \overline{\{ u\}}    \big) \equiv  \mathscr{E}_7 \big(       e ,   \overline{\{ u\}}         \big) \equiv   \big\{ e \in \mathcal{E}  \big(   \overline{\{ u \}}   \big) :      e \in  \mathcal{E} \big(   u \cap \Gamma    \big)  \cap \mathcal{E} \big( \overline{\{ u\}} \big)      \big\}         \text{, } 
       \end{align*}

       \noindent where,
       
       \begin{align*}
 \mathcal{E}    \supsetneq \mathcal{E} \big(   u \cap \Gamma    \big)   \equiv \big\{     e \in \mathcal{E} :  e \cap u \cap \Gamma \neq \emptyset                       \big\} \text{, } 
              \end{align*}
       \noindent which, from previous observation, 
       
       \begin{align*}
        \mathscr{E}^{7}_{u , \Gamma , \geq h } \equiv \mathscr{E}_{u , \Gamma , \geq h}   \equiv \big\{              e \in \mathcal{E} :                       \big\{  e \cap \big\{ u \overset{\mathscr{Z} , \geq h}{\longleftrightarrow} \Gamma \big\}  \neq \emptyset \big\} \cup \big\{ d_{u , \Gamma}  \big\}                                          \big\}               \text{, }  
       \end{align*}
       
       \noindent 
       
       \noindent for the smallest length of the path realized by the same infimum $ d_{u , \Gamma}$. Concluding,
       
       \begin{align*}
             \textbf{P}_h \big[   e \in \mathscr{E}_{u , \Gamma , \geq h}                           \big] \overset{\mathrm{'edge\text{ }  opening'}}{\Longleftrightarrow}       \textbf{P}_h \big[       e \in \mathscr{E}_{u , \Gamma , < h}                       \big]   \text{ } \text{ . }
       \end{align*}

          \noindent As a result,
       
       \begin{align*}
      \big| \omega \cap f_h \big( \omega \big) \big| < \mathscr{D}^{\prime\prime}
      \text{, } 
       \end{align*}

       \noindent for suitable, strictly positive $\mathscr{D}^{\prime\prime}$. Hence $\overline{\{ u\}}$ can be reconstructed from $f_h \big( \omega \big)$ as previously described.
       
            \noindent \underline{Case D}

            \begin{itemize}
            
            \item[$\bullet$] \underline{Modification one} (\textit{closing edges of the metric-graph GFF configuration}). For the final case, we refer to \underline{Modification one} for \underline{Case C} included in the case above first. With the same objects introduced for \underline{Modification one}, also introduce the following two objects, the first of which is given by,
            
            \begin{align*}
            \mathscr{E}^7_u \equiv   \mathscr{E}^7 \big( e , \Gamma ,  \pi_u \big)   \equiv \mathscr{E}^7 \big( e ,  \pi_u \big)  \equiv  \big\{  e \in \mathcal{E}       :         e \in \mathcal{E} \big(  \pi_u  \cap    \Gamma  \big)   \cap     \mathcal{E} \big( u \big)            \big\}        \text{ } \text{ , }
            \end{align*}

            \noindent and the second of which is similarly given by the edge-set,
            
             \begin{align*}  \mathscr{E}^7_v \equiv   \mathscr{E}^7 \big( e , \Gamma ,  \pi_v \big)   \equiv \mathscr{E}^7 \big( e ,  \pi_v \big)  \equiv  \big\{  e \in \mathcal{E}       :         e \in \mathcal{E} \big(  \pi_v  \cap    \Gamma  \big)   \cap     \mathcal{E} \big( v \big)            \big\}   
                \text{ } \text{ , }
            \end{align*}
            
            \noindent implying,

            \begin{align*}
               \textbf{P}_h \big[  e \in    \mathscr{E}^7_u  \cap \mathscr{E}^7_v  \cap \mathscr{E}^7  : \big\{         e               \overset{\mathscr{Z} , < h} {\longleftrightarrow}          \overline{\{ u\}}                  \big\}     \big]   > 0           \overset{\mathrm{'edge\text{ }  closing'}}{\Longleftrightarrow}  \textbf{P}_h \big[  e \in    \mathscr{E}^7_u  \cap \mathscr{E}^7_v  \cap \mathscr{E}^7  :           \big\{                      e   \overset{\mathscr{Z} , \geq h} {\longleftrightarrow}      \overline{\{ u\}}          \big\}             \big] > 0            \text{, } 
            \end{align*}
            
           \noindent with the intersection of three edge sets, $\mathscr{E}^7_u \cap \mathscr{E}^7_u \cap \mathscr{E}^7$,
           
                       \begin{align*}    \mathscr{E}^7_u \cap \mathscr{E}^7_u \cap \mathscr{E}^7 \equiv \big\{   e \in \mathcal{E} \big( \mathcal{S} \big) :  e \in \mathcal{E} \big( \pi_u \cap \Gamma \big) \cap \mathcal{E} \big( \pi_v \cap \Gamma \big) \cap    \mathcal{E} \big( u \big) \cap \mathcal{E} \big( v \big)           \big\}       
                \text{. }
            \end{align*}

            \item[$\bullet$] \underline{Modification two} (\textit{opening edges from the shortest path of the metric-graph GFF configuration after opening edges}). From \underline{Modification two} introduced in \underline{Case C}, introduce the following modification to the edge set,
            
            \begin{align*}
              \mathscr{E}^8_u            \equiv     \mathscr{E}^7_u          \text{, } 
            \end{align*}
            
            \noindent in addition to the definition of the following path,

            \begin{align*}
             \rho_{u, \mathcal{S}} \equiv \rho_u \equiv \underset{\Gamma \subset \mathcal{S}}{\underset{u \in \mathcal{V} \big(\mathcal{S}\big)}{\mathrm{inf}}}  \big\{    \rho_{\mathcal{S}} \big(  u              ,    \Gamma          \big)      \big\}  \text{ } \text{ , }
            \end{align*}
            
            \noindent which for the shortest such path, implies that \underline{Modification two} takes the form, for path realizations $\gamma$ of the shortest path between $u$ and $\Gamma$,
            
            \begin{align*}
             \textbf{P}_h \big[   e \in \mathcal{E}^8_u : \big\{  \gamma \equiv \rho_u   \big\} \cup      \big\{   u  \overset{\mathscr{Z} , \geq h}{\longleftrightarrow}     \Gamma     \big\}              \big]  \overset{\mathrm{'edge\text{ }  opening'}}{\Longleftrightarrow}   \textbf{P}_h \big[              e \in \mathcal{E}^8_u : \big\{  \gamma \equiv \rho_u   \big\} \cup        \big\{  u   \overset{\mathscr{Z} , < h}{\longleftrightarrow}     \Gamma     \big\}           \big]  \text{, } 
            \end{align*}
            
            \noindent as opening the edges for the metric-graph GFF entails the change to the level-set crossing event,

            \begin{align*}
                  \big\{  u   \overset{\mathscr{Z} , < h}{\longleftrightarrow}     \Gamma     \big\}              \text{, } 
            \end{align*}
            
            \noindent a modification to the height of the metric-graph GFF crossing event.
            
            \item[$\bullet$] \underline{Modification three} (\textit{opening edges of the metric-graph GFF configuration once more}). Introducing the following modification to edge sets used in \underline{Modification two} of \underline{Case D} implies,

              \begin{align*}
             \textbf{P}_h \big[   e \in \mathcal{E}^8_u : \big\{  \gamma \equiv \rho^{\prime}_v   \big\} \cup      \big\{   u  \overset{\mathscr{Z} , \geq h}{\longleftrightarrow}     \Gamma     \big\}              \big]  \overset{\mathrm{'edge\text{ }  opening'}}{\Longleftrightarrow}   \textbf{P}_h \big[              e \in \mathcal{E}^8_u : \big\{  \gamma \equiv \rho^{\prime}_v   \big\} \cup        \big\{  u   \overset{\mathscr{Z} , < h}{\longleftrightarrow}     \Gamma     \big\}           \big]  \text{, } 
            \end{align*}
            
            \noindent instead for the distance minimizing path $\gamma$ between $v$ and $\Gamma \cup \rho$, as,
            
              \begin{align*}
             \rho^{\prime}_{v, \mathcal{S}} \equiv \rho^{\prime}_v \equiv \underset{\Gamma \cup \rho  \subset \mathcal{S}}{\underset{v \in \mathcal{V} \big(\mathcal{S}\big)}{\mathrm{inf}}}  \big\{    \rho_{\mathcal{S}} \big(  v              ,    \Gamma  \cup \rho        \big)      \big\}  \text{ } \text{ . }
            \end{align*}

            \end{itemize}

              \end{itemize}
       
   \noindent As a result,
       
       \begin{align*}
      \big| \omega \cap f_h \big( \omega \big) \big| < \mathscr{D}^{\prime\prime\prime}
      \text{, } 
       \end{align*}

       \noindent for suitable, strictly positive $\mathscr{D}^{\prime\prime\prime}$. Hence $\overline{\{ u\}}$ can be reconstructed from $f_h \big( \omega \big)$ as previously described.

       \bigskip

       \noindent This concludes the arguments to demonstrate that the upper bound displayed in ($\bar{\Gamma}$ \textit{Probability II}) holds, from which we conclude the argument for establishing quasi-multiplicativity. \boxed{}
       
       \section{The Villain model, and its cable system counterpart, akin to the metric-graph extension of the GFF}
       
       \noindent We now turn towards defining another model of interest, in which the authors of {\color{blue}[12]} propose that it can be of interest co construct IIC-type limits also for the Villain model, from the expectation value, for $k \equiv 1$, $k \equiv 2$, and some $\theta >0$,
       
       \begin{align*}
      \langle       \mathrm{exp} \big(   i k  \theta_0    \big)      \rangle_{\mathrm{IIC}_k \propto\textbf{P}^{\mathscr{V}} }     \text{, } 
       \end{align*}

       \noindent of an observable, for $k \equiv 1$, with respect to the the incipient infinite cluster measure that is proportional to the base measure $\textbf{P}^{\mathscr{V}} \big[ \cdot \big]$ for the Villain model which will be introduced formally below. Denote the finite subgraph of $\textbf{Z}^d$ with connected $\mathcal{G} \equiv G \equiv \big( V , E \big)$. As with the metric-graph extension for the GFF in three dimensions presented at the beginning of the paper, denote the metric-graph from $\mathcal{G}$ with $\widetilde{\mathcal{G}} \equiv \widetilde{G} \equiv  \big( \widetilde{V} , \widetilde{E} \big)$.

  \subsection{Villain objects}

       \noindent \textbf{Definition} \textit{14} (\textit{Villain model Hamiltonian}). Introduce,

       \begin{align*}
    \mathcal{H} \big( \theta \big) \equiv  \underset{n \in \textbf{Z}}{\prod}    \big( - \beta \big( \theta_x - \theta_y + 2 \pi n \big)^2 \big)      \text{, } 
       \end{align*}
       
       \noindent corresponding to the Hamiltonian of the Villain model, for the following configurations drawn from spin space $\textbf{U}$,

       \begin{align*}
     \big( u_x \big)_{x \in V}   \equiv \big( \mathrm{exp} \big( i \theta_x \big) \big)_{x \in V} \in \Omega^{\mathrm{Villain}} \simeq \frac{\textbf{R}}{2 \pi \textbf{Z}} \simeq \textbf{U}  \text{, } 
       \end{align*}
       
       \noindent at inverse temperature $\beta>0$, for $\theta_x , \theta_y > 0$.

       \bigskip
       
       \noindent \textbf{Definition} \textit{15} (\textit{Villain model probability measure from the Hamiltonian}). Introduce,
       
       \begin{align*}
      \big(  \textbf{P}^{\mathscr{V}}_{\mathcal{G}}  \big[    \omega   \big] \big)^{\theta}  \equiv  \textbf{P}^{\mathscr{V}}_{\mathcal{G}}  \big[    \omega   \big]  \equiv  \mathcal{Z}^{-1} \mathrm{exp} \big(  - \mathcal{H} \big( \theta \big)     \big) \underset{x \in V}{\prod} \mathrm{d} \theta_x \equiv   \mathcal{Z}^{-1} \underset{\{ xy\} \in E}{\prod} \mathrm{exp} \big\{           -  \underset{n \in \textbf{Z}}{\prod}    \big( - \beta \big( \theta_x - \theta_y + 2 \pi n \big)^2 \big)        \big\} \underset{x \in V}{\prod} \mathrm{d} \theta_x  \\ \overset{\beta \longrightarrow + \infty }{\approx}    \mathcal{Z}^{-1} \underset{\{ xy\} \in E}{\prod}  \mathrm{exp} \big\{           -  \underset{n \in \textbf{Z}}{\prod}    \big( - \frac{\beta}{2} \big( \mathrm{cos} \big( \theta_x - \theta_y \big) \big)        \big\} \underset{x \in V}{\prod} \mathrm{d} \theta_x       \text{, } 
       \end{align*}

       \noindent corresponding to the probability of the Villain model, with configurations from sample space $\Omega^{\mathrm{Villain}}$ provided in the previous definition, with partition function,

       \begin{align*}
    \mathcal{Z}^{\mathscr{V}}_{\mathcal{G}} \big( \theta \equiv 1,  \theta_x , \theta_y , \beta \big)   \equiv   \mathcal{Z}^{\mathscr{V}}_{\mathcal{G}} \big( \theta_x , \theta_y , \beta \big)  \equiv \mathcal{Z} \big( \theta_x , \theta_y , \beta \big)  \equiv \mathcal{Z} \equiv    \underset{x , y \in \mathcal{G}}{\sum } \underset{\{ xy\} \in E}{\prod} \mathrm{exp} \big\{           -  \underset{n \in \textbf{Z}}{\prod}    \big( - \beta \big( \theta_x - \theta_y + 2 \pi n \big)^2 \big)        \big\}  \\ \times   \underset{x \in V}{\prod} \mathrm{d} \theta_x  \text{, } 
       \end{align*}
       
       \noindent for boundary conditions $\theta \equiv 1$.
       
       \bigskip
       
       \noindent \textbf{Definition} \textit{16} (\textit{transition probability for Brownian motion over the Villain model sample space for the cable system representation}). Introduce an exponential with a similar Hamiltonian to that provided in \textbf{Definition} \textit{14}, namely,

       \begin{align*}
              p_t \big( \theta_1 , \theta_2 \big) \equiv    \frac{1}{\sqrt{2 \pi t}} \underset{n \in \textbf{Z}}{\sum}    \mathrm{exp} \big\{     - \big( 2 t \big)^{-1} \big( \theta_1 - \theta_2 + 2 \pi n \big)^2 \big)        \big\}     \text{, } 
       \end{align*}

       \noindent corresponding to the \textit{probability transition function} for Brownian motion over $\Omega^{\mathrm{Villain}}$, for $t<1$, or $t>1$.  
       
       \bigskip
       
       \noindent \textbf{Definition} \textit{17} (\textit{cable system representation of the Villain model Hamiltonian}). Introduce,
       
       \begin{align*}
      \mathcal{H}^{\mathrm{C}-\mathrm{S}} \big( \theta \big) \equiv   \mathcal{H} \big( \theta \big) \equiv   \mathrm{ln} \big\{   \big( 2 \pi t \big)^{\frac{|V|}{2}}        \underset{\{ x y \} \in \widetilde{E}}{\prod}   p_t \big( \theta_1 , \theta_2 \big)      \big\}     \text{, } 
       \end{align*}

       \noindent corresponding to the representation of the Hamiltonian for the Villain model under the cable system.

       \bigskip
       
       \noindent \textbf{Definition} \textit{18} (\textit{cable system probability measure for the extended Villain model}). Introduce,
      
      \begin{align*}
          \big( \widetilde{\textbf{P}^{\mathscr{V}}_{\mathcal{\widetilde{G}}}} \big[       \widetilde{\omega}        \big]  \big)^{\widetilde{\theta}} \equiv   \widetilde{\textbf{P}^{\mathscr{V}}_{\mathcal{\widetilde{G}}}} \big[       \widetilde{\omega}        \big]   \equiv \widetilde{\mathcal{Z}^{-1}} \mathrm{exp} \big\{  - \mathcal{H} \big( \theta \big)      \big\}   \equiv \widetilde{\mathcal{Z}^{-1}} \big( 2 \pi t \big)^{\frac{|V|}{2}} \underset{\{ xy \} \in E}{\prod}   p_t \big( \theta_x , \theta_y \big)         \text{, } 
      \end{align*}

       \noindent corresponding to the probability measure for the extended Villain model, with partition function,

       \begin{align*}
   \widetilde{\mathcal{Z}^{\mathscr{V}}_{\widetilde{\mathcal{G}}}}  \big( \widetilde{\theta} \equiv 1, \theta_1 , \theta_2 , t  \big)    \equiv   \widetilde{\mathcal{Z}^{\mathscr{V}}_{\widetilde{\mathcal{G}}}}  \big( \theta_1 , \theta_2 , t  \big) \equiv    \widetilde{\mathcal{Z}^{\mathscr{V}}_{\widetilde{\mathcal{G}}}}  \equiv \widetilde{\mathcal{Z}} \equiv      \underset{x , y \in \widetilde{\mathcal{G}}}{ \sum} \big( 2 \pi t \big)^{\frac{|V|}{2}} \underset{\{ xy \} \in E}{\prod}   p_t \big( \theta_x , \theta_y \big)             \text{, } 
       \end{align*}
       
       \noindent for boundary conditions $\widetilde{\theta} \equiv 1$.
       
       \bigskip
       
       \noindent From these objects, we proceed to also introduce connectivity events for the extended Villain model.

       \bigskip
       
       \noindent \textbf{Definition} \textit{19} (\textit{connectivity events}). Introduce, for $S \subset \textbf{U}$, and $\theta \in C \big( \widetilde{\mathcal{G}} , \textbf{U} \big)$,
       
       \begin{align*}
           \big\{  x           \overset{S}{\longleftrightarrow}      y       \big\}             \text{, } 
       \end{align*}

       \noindent corresponding to a connectivity event, between $x,y\in \widetilde{\mathcal{G}}$ along some path $\gamma$, so that $\theta \big( \gamma \big) \subset S$.

       \bigskip
       
       \noindent Additionally, besides the Villain crossing event denoted above, we also introduce \textit{Villain clusters} for the extended model.

       \bigskip

       \noindent \textbf{Definition} \textit{20} (\textit{extended Villain clusters}). Introduce, for $\Lambda \subsetneq \textbf{Z}^d$ with boundary $\partial \Lambda$,
       
       \begin{align*}
       \mathcal{C}^{\mathrm{Villain}} \equiv \mathcal{C} \equiv \big\{    y \in \widetilde{\mathcal{G}} :    y \overset{S}{\longleftrightarrow}   \partial \Lambda  \big\}     \text{, } 
       \end{align*}
       
       \noindent corresponding to clusters in the extended Villain model, within $S$.
       
       \bigskip
       
       \noindent \textbf{Definition} \textit{21} (\textit{linear heat kernel in the extended Villain model with absorbing boundary}). Introduce,
       
       \begin{align*}
                p^{[-\frac{\pi}{2} , \frac{\pi}{2}]}_t \big( \theta_1 , \theta_2    \big)    \equiv    \frac{1}{\sqrt{2 \pi t}}    \underset{ n \in \textbf{Z}}{\sum  } \bigg(      \mathrm{exp} \bigg(    - \big( 2 t \big)^{-1}  \big( \theta_1 - \theta_2 + 2 \pi n       \big)^2 \bigg)        -  \mathrm{exp} \bigg(         - \big( 2 t \big)^{-1}   \big(   \theta_1 + \theta_2 + \big( 2 n - 1 \big) \pi      \big)^2     \bigg)    \bigg)    \text{, } 
       \end{align*}

       \noindent corresponding to the heat kernel of the extended Villain model, obtained from repeated reflection principle as discussed in {\color{blue}[12]}.
       
       \bigskip
       
       \noindent From the linear heat kernel, also introduce the following result.
       
       \bigskip

       \noindent \textbf{Definition} \textit{22} (\textit{extended connectivity event probability}). Given $x,y \in \widetilde{\mathcal{G}}$, and $z \in \{ x y \}$, as introduced in \textbf{Definition} \textit{19}, introduce,
       
       \begin{align*}
        \big\{   \mathrm{exp} \big( i \theta_x \big)   \equiv   u \big( x \big)   \overset{\{ \pm i \}^c }{\longleftrightarrow}         u \big( y \big)    \equiv  \mathrm{exp} \big( i \theta_y \big)      \big\}               \text{, } 
       \end{align*}

       \noindent corresponding to the \textit{connectivity event} between $u(x)$ and $u(y)$ in $\{ \pm i \}^c$.

       \bigskip

       \noindent Finally, introduce the following $\textbf{Lemma}$, from {\color{blue}[12]}, for relating the probability of the connectivity event above, namely,
       
       \begin{align*}
            \widetilde{\textbf{P}^{\mathscr{V}}_{\mathcal{\widetilde{G}}}} \big[   u \big( x \big)   \overset{\{ \pm i \}^c }{\longleftrightarrow}         u \big( y \big)            \big]        \text{, } 
       \end{align*}

       \noindent with the linear heat kernel.

  \bigskip
  
  \noindent \textbf{Lemma} \textit{3.1} (\textit{equating the probability of an extended connectivity event occurring with a normalized linear heat kernel}, ({\color{blue}[12]}, \textbf{Lemma} \textit{3.1})). Given parameters specified in \textbf{Definition} \textit{22}, the probability of an \textit{extended connectivity event} occurring between $u(x)$ and $u(y)$ is,
  
  \begin{align*}
     0 \text{ if } \mathrm{cos} \big( \theta_x \big) \mathrm{cos} \big( \theta_y \big) \leq 0 \\   \frac{p^{[-\frac{\pi}{2} , \frac{\pi}{2}]}_t \big( \theta_x , \theta_y    \big)}{p_t \big( \theta_x , \theta_y \big) }    \text{ } \text{ if } \text{ }   \theta_x , \theta_y \in \big( - \frac{\pi}{2} , \frac{\pi}{2} \big) \text{, }  \\               \frac{p^{[-\frac{\pi}{2} , \frac{\pi}{2}]}_t \big( \theta_x - \pi  , \theta_y - \pi    \big)}{p_t \big( \theta_x , \theta_y \big) }  \text{ } \text{ if } \text{ }       \theta_x , \theta_y \in \big( \frac{\pi}{2} ,  \frac{3\pi}{2} \big) \text{ }  \mathrm{mod} \big( 2 \pi \big)   \text{, }   
  \end{align*}
  
  \noindent depending upon the conditions enforced on $\theta_x$ and $\theta_y$ from each of the three possible values of the \textit{extended connectivity event} occurring.
  
    \bigskip

\noindent As with the GFF which satisfies (FKG), the Villain model, as a variant of the XY model, also satisfies (FKG), the statement of which is provided below.

\bigskip

\noindent \textbf{Proposition} \textit{Villain} (\textit{FKG property of the Villain model from FKG property of the XY model}, ({\color{blue}[33]}, \textbf{Proposition} \textit{A.1}). The XY probability measures $\rho \big[ \cdot \big]$ are FKG.

    \bigskip
    
    \noindent From the probability of an \textit{extended connectivity event occurring} from the linear heat kernel, the conditional probability that $\big\{ x \overset{\{ \pm i\}^c}{\longleftrightarrow} y \big\}$ occurs also admits a particular representation which is given below.

       \bigskip
       
       \noindent \textbf{Proposition} \textit{Villain-1} (\textit{conditional probability of a connectivity event occurring}, {\color{blue}[12]}). Conditionally upon $\theta|_{\mathcal{G}}$, the probability,

       \begin{align*}
          \widetilde{\textbf{P}^{\mathscr{V}}_{\mathcal{\widetilde{G}}}} \big[   u \big( x \big)   \overset{\{ \pm i \}^c }{\longleftrightarrow}         u \big( y \big)  \text{ }   |     \text{ }  \theta|_{\mathcal{G}}          \big]                 \text{, } 
       \end{align*}
       
       \noindent is equivalent to,

       \begin{align*}
          \mathrm{max} \big( 0 , p + q - 1 \big) \leq c \leq \mathrm{min} \big( p , q \big)                    \text{, } 
       \end{align*}

       \noindent for parameters $c,p,q>0$.
       
       \bigskip
       
       \noindent Before presenting further arguments, under previously stated assumptions, 
       
       \begin{align*}
       \big(  \textbf{P}^{\mathscr{V}}_{\mathcal{G}}  \big[ \big( u_x \big)_{x \in V}  \big] \big)^{\theta}  \equiv    \big(  \textbf{P}^{\mathscr{V}}_{\mathcal{G}}  \big[ \big( \overline{u}_x \big)_{x \in V}  \big] \big)^{\theta}     \text{, } 
       \end{align*}

       \noindent for boundary conditions $\theta$.
       
       \bigskip
       
       \noindent Finally, it is also worth mentioning that, as is the case for the metric-graph GFF, that the Villain model also satisfies a version of the strong Markov property (see \textbf{Lemma} \textit{2.1}, from {\color{blue}[12]}).

       \subsection{Formulating an IIC-type limit in the Villain context}

       \noindent From prerequisites in the previous section, we also introduce the following two items.
       
       \bigskip

       \noindent \textbf{Proposition} \textit{Villain-2} ({\color{blue}[12]}, \textbf{Proposition} \textit{3.2}). For each $t_0>0$, there exists $\delta \equiv \delta \big( t_0 , \mathrm{deg}(x) \big) >0$, such that, if $t \big( xy \big) \geq t_0$ for all neighbors $y \sim x$, then,
       
       \begin{align*}
         \delta \leq \frac{\langle \mathrm{cos} \big( \theta_x \big) \rangle}{\textbf{P} \big[  x         \overset{\{ \pm i \}^c }{\longleftrightarrow} \partial       \big] } \leq 1  \text{, } 
       \end{align*}
       
       \noindent for $\partial \equiv \partial V$, and $\textbf{P}_{\mathscr{G}} \big[ \cdot \big]  \equiv \textbf{P} \big[ \cdot \big]$.
       
       \bigskip
      
         \noindent \textbf{Proposition} \textit{Villain-3} ({\color{blue}[12]}, \textbf{Proposition} \textit{3.3}). For each $t_0>0$, there exists $\delta \equiv \delta \big( t_0 , \mathrm{deg}(x) \big) >0$, such that, if $t \big( xy \big) \geq t_0$ for all neighbors $y \sim x$, then, 
         
       \begin{align*}
        \delta \leq \frac{\langle \mathrm{cos} \big( 2 \theta_x \big) \rangle}{\textbf{P} \big[  x         \overset{\{ \pm \xi \}^c }{\longleftrightarrow} \partial     , x         \overset{\{ \pm \bar{\xi} \}^c }{\longleftrightarrow} \partial    \big] } \leq 1     \text{, } 
       \end{align*}
       
       \bigskip
       
         \noindent for $\partial \equiv \partial V$, and $\textbf{P}_{\mathscr{G}} \big[ \cdot \big]  \equiv \textbf{P} \big[ \cdot \big]$.

         \bigskip
         
         \noindent From each \textbf{Proposition} above, introduce the two limits for the Villain model IIC-type limit which will be shown to coincide.

         \bigskip

         \noindent \textbf{Definition} \textit{23} (\textit{cylinder events for the Villain model}). Akin to \textit{cylinder events} for Bernoulli percolation in two-dimensions, and for the metric-graph GFF in three-dimensions, a \textit{cylinder event} for the Villain model is dependent upon a finite number of vertices over $\textbf{Z}^2$, in addition to the $\theta$ interaction, as,

         \begin{align*}
          \mathscr{E} \big( \theta \big) \equiv \mathscr{E}   \text{, } 
         \end{align*}

         \noindent for $\theta \in [ 0 , 2 \pi ]$.
         
         \bigskip

         \noindent \textbf{Definition} \textit{24} (\textit{Villain model equivalent of the conditional crossing limit for Bernoulli at criticality from \textbf{Definition} \textit{6}}). Introduce, for $\theta$ interactions coinciding with $\{ \pm i \}^c$, and $w \in V$, the probability,

         \begin{align*}
        \mathcal{P}^{\prime}_1 \equiv    \textbf{P} \big[  x         \overset{\{ \pm i \}^c }{\longleftrightarrow} \partial   \big|    | \mathcal{C} \big(  w             \big)  | = + \infty        \big]  \equiv    \textbf{P} \big[ \mathscr{E} \big( \pm \frac{\pi}{2} \big)   \big|    | \mathcal{C} \big(  w             \big)  | = + \infty         \big]        \text{, } 
         \end{align*}
         
         \noindent corresponding to the first limiting ratio of the IIC-type limit for the Villain model.

         \bigskip
         
         \noindent \textbf{Definition} \textit{25} (\textit{Villain model equivalent of the limit of conditionally defined crossing probabilities for Bernoulli percolation at supercriticality from \textbf{Definition} 7}). Introduce, for $\theta$ interactions coinciding with $\{ \pm \xi\}^c$, and $\{ \pm \bar{\xi} \}^c$, as $\xi , \bar{\xi} \searrow i $, and $w \in V$, the intersection probability,

         \begin{align*}
\mathcal{P}^{\prime\prime\prime}_2 \equiv  \underset{\xi , \bar{\xi} \searrow i}{\mathrm{lim}}  \mathcal{P}^{\prime\prime}_2 \equiv   \underset{\xi , \bar{\xi} \searrow i}{\mathrm{lim}}  \textbf{P} \big[    x         \overset{\{ \pm \xi \}^c }{\longleftrightarrow} \partial     , x         \overset{\{ \pm \bar{\xi} \}^c }{\longleftrightarrow} \partial  \big|    | \mathcal{C} \big(  w             \big)  | = + \infty                                \big]    \equiv  \underset{\xi , \bar{\xi} \searrow i}{\mathrm{lim}} \textbf{P} \big[  \mathscr{E} \big( \pm i^{-1} \mathrm{ln} \big\{       \xi           \big\}  \big)  ,  \mathscr{E} \big(\pm  i^{-1} \mathrm{ln} \big\{        \bar{\xi}         \big\} \big) \big|  | \mathcal{C} \big(  w             \big)  | = + \infty                  \big]   \text{ , } 
         \end{align*}

         \noindent corresponding to the second limiting ratio of the IIC-type limit for the Villain model.

         \bigskip

         \noindent From \textbf{Definition} \textit{23} and \textbf{Definition} \textit{24}, the next item below asserts the existence of the IIC-type limit for the Villain model.

         \bigskip
         
         \noindent \textbf{Theorem} \textit{7} (\textit{establishing equivalence between  the two limiting ratio of the Villain IIC-type limit, analogous to \textbf{Theorem} 3 for Bernoulli percolation, and to \textbf{Theorem} 5 for the metric-graph GFF}). Given the assumptions provided above, the two limits, respectively from \textbf{Definition} \textit{23}, and \textit{24}, satisfy,

         \begin{align*}
              \mathcal{P}^{\prime}_1 \equiv \mathcal{P}^{\prime\prime\prime}_2          \text{, } 
         \end{align*}
         
         \noindent which, by Kolmogorov extension, implies,

         \begin{align*}
            \nu_{\theta}   \big[         \forall w \in V , \exists ! \mathcal{C} \big( w \big) : | \mathcal{C} \big( w \big) | = + \infty                     \big] \equiv  1   \text{, } 
         \end{align*}
         
         \noindent upon denoting the common value between $\mathcal{P}^{\prime}_1$ and $\mathcal{P}^{\prime\prime\prime}_2$ with,
         
         \begin{align*}
         \nu_{\theta} \big[ \mathscr{E} \big]              \text{. }
         \end{align*}

         \noindent To prove \textbf{Theorem} \textit{7} pertaining to the Villain IIC-type limits, we appeal to arguments provided in the streamlined presentation of Kesten's seminal arguments, from {\color{blue}[24]}, in {\color{blue}[3]}.

         \bigskip
         
         \noindent \textit{Proof of Theorem 7}. From arguments originally introduced for Bernoulli percolation in two-dimensions, which were adapted in \textbf{Theorem} \textit{5} for the three-dimensional metric-graph GFF, we introduce the following, namely,
         
         \begin{align*}
             \textbf{P} \big[       \mathscr{E} \text{ } | \text{ }            x \overset{\{ \pm \xi \}^c }{\longleftrightarrow}   \partial                                 \big]          \longrightarrow  \nu_{\theta} \big[ \mathscr{E} \big]  \text{, } 
         \end{align*}

         \noindent corresponding to the Kolmogorov extension result provided in $(\textit{1})$. Corresponding to $(\textit{2})$, introducing an analogous quantity yields,

       \begin{align*}
         \epsilon_i \equiv   \underset{\mathscr{F}_1 , \mathscr{F}_2 \in \mathcal{A}_i}{\underset{\xi \in [             0 , i       ]}{\mathrm{sup}}} \text{ } \textbf{P} \big[          \big\{         \mathscr{F}_1      \underset{\mathcal{A}_i}{\overset{\{ \pm \xi \}^c }{\longleftrightarrow} }      \mathscr{F}_2 \big\}^c    \text{ } | \text{ } \partial_i \overset{\{ \pm \xi \}^c }{\longleftrightarrow} 
    \partial_{i+1}    \big]                             \text{, } 
       \end{align*}

       \noindent for,

       \begin{align*}
      \partial_{i} \equiv \partial  \mathcal{B}_i        \text{, } \\ \partial_{i+1} \equiv  \partial  \mathcal{B}_{i+1}    \text{, } 
 \end{align*}
         
       \noindent and $\mathscr{F}_1 , \mathscr{F}_1 \subsetneq  F \big(   \mathcal{A}_i  \big) \subsetneq F \big( \textbf{Z}^2 \big)$. From the two quantities introduced above, the remainder of the arguments from \textbf{Theorem} \textit{5} can be executed to demonstrate that,

       \begin{align*}
      \mathcal{P}^{\prime}_1 \equiv \textbf{P} \big[ \mathscr{E} \big( \pm \frac{\pi}{2} \big)   \big|    | \mathcal{C} \big(  w             \big)  | = + \infty         \big]     \equiv      \underset{\xi , \bar{\xi} \searrow i}{\mathrm{lim}} \textbf{P} \big[  \mathscr{E} \big( \pm i^{-1} \mathrm{ln} \big\{       \xi           \big\}  \big)  ,  \mathscr{E} \big(\pm  i^{-1} \mathrm{ln} \big\{        \bar{\xi}         \big\} \big) \big|  | \mathcal{C} \big(  w             \big)  | = + \infty                  \big]    \equiv      \mathcal{P}^{\prime\prime\prime}_2      \text{, } 
       \end{align*}

       \noindent upon following the arguments from $(1^{*})$ to $(2^{*})$ for obtaining an upper bound for an intersection of a Villain crossing probability with a cylinder event due to quasi-multiplicativity. Following a similar lower bound obtained in $(2^{*})$, expressing a conditionally defined crossing event for the Villain model, as provided in $(3^{*})$ for the three-dimensional metric-graph GFF, allows for one to introduce another counterpart for the quantity provided in $(4^{*})$ for the Villain model. Proceeding from $(4^{*})$, introducing similarly defined quantities allows for one to, from a crossing event for the Villain model analogous to that introduced for the three-dimensional metric-graph GFF in $(5^{*})$, an exponential upper bound similar to $(7^{*})$ following an analogous lower bound to $(6^{*})$, from which we conclude the argument. \boxed{}

    \section{Quasi-multiplicativity of crossing probabilities for the two-dimensional Villain IIC-type limit}

    \noindent The arguments for establishing quasi-multiplicativity over slabs for crossing probabilities in Bernoulli peroclation provided in {\color{blue}[3]} directly carry over to the Villain model, and are hence excluded (refer to \textbf{Theorem} \textit{3.1} and \textbf{Lemma} \textit{3.2} provided in {\color{blue}[3]} for details).

       \section{Acknowledgments}

       \noindent The author would like to thank Philippe Sosoe for comments.
       
\section{References}

\noindent [1] Abacherli, A. $\&$ Sznitman, A-S. Level-set percolation for the Gaussian free field on a transient tree. Annales de I'Institut Henri Poincare - Probabilites et Statistiques \textbf{54}(1), 173-201 (2018).

\bigskip

\noindent [2] Abacherli, A. $\&$ Cernyi, J. Level-set percolation of the Gaussian free field on regular graphs II: finite expanders. Electronic Journal Probability \textbf{125}(130) 1-39 (2020). 

\bigskip

\noindent [3] Basu, D., Sapozhnikov, A. Kesten's incipient infinite cluster and quasi-multiplicativity of crossing probabilities. \textit{Elec Comm Commun. Probab., 22:Paper No. 26, 12} (2017).

\bigskip

\noindent [4] Basu, D., Sapozhnikov, A. Crossing probabilities for critical Bernoulli percolation on slabs. \textit{arXiv 1512.05178, v2}.

\bigskip

\noindent [5] Berestycki, N., Powell, E. $\&$ Ray, G. A characterisation of the Gaussian free field. Probability Theory and Related Fields \textbf{176}, 1259-1301 (2020).

\bigskip

\noindent [6] Biskup, M. Extrema of discrete two-dimensional Gaussian free field. arXiv:1712.09972v4 (2019).

\bigskip

\noindent [7] Borodin, A. $\&$ Gorin, V. General $\beta$- Jacobi Corners Process and the Gaussian Free Field. Communications on Pure and Applied Mathematics \textbf{LXVIII} 1774-1844 (2015).

\bigskip

\noindent [8] Bramson, W. Zeitouni, O. Tightness of the recentered maximum of the discrete Gaussian free field. Communications on Pure and Applied Mathematics \textbf{LXV} 0001-0020 (2012).

\bigskip

\noindent [9] Cao, X. $\&$ Santachiara, R. Level-set percolation in two-dimensional Gaussian free field. arXiv:2012.09570v3 (2021).

\bigskip

\noindent [10] Ding, J. $\&$ Li, L. Chemical distances for percolation of planar Gaussian free field and critical random walk loop soups. Communications in mathematical physics, \textbf{360} 523-553 (2018).

\bigskip

\noindent [11] Ding, J. $\&$ Wirth, M. Percolation for level-sets of the Gaussian free field on metric graphs. arXiv: 1807.11117v2 (2019).

\bigskip

\noindent [12] Dubedat, J., Falconet, H. Random clusters in the Villain and XY models. \textit{arXiv:2210.03620, v1} (2022).

\bigskip

\noindent [13] Duminil-Copin, H., Goswami, S., Rodriguez, P-F., $\&$ Severo, F. Equality of critical parameters for percolation of Gaussian free field level-sets. arXiv: 2002.077735v1 (2020).

\bigskip

\noindent [14] Duminil-Copin, H., Goswami, S., Raoufi, A., Severo, F. $\&$ Yadin, A. Existence of phase transition for percolation using the Gaussian free field. arXiv: 1816.07733v2 (2020).

\bigskip

\noindent [15] Drewitz, A. $\&$ Rodriguez, P-F. High-dimensional asymptotics for percolation of Gaussian free field level sets. arXiv: 1310.1041v3 (2015).

\bigskip

\noindent [16] Drewitz, A., Prevost, A., $\&$ Rodriguez, P-F. Critical exponents for a percolation model on transient graphs. arXiv: 2101.05801v2 (2021).

\bigskip

\noindent [17] Drewitz, A., Prevost, A. $\&$ Rodriguez, P-F. The sign clusters of the massless Gaussian free field percolate on $\textbf{Z}^d$, $d \geq 3$ (and more). arXiv: 1708.03285v2 (2018).

\bigskip

\noindent [18] Drewitz, A., Prevost, A. $\&$ Rodriguez, P-F. Geometry of Gaussian free field sign clusters and random interlacements. arXiv: 1811.05970v1 (2018).

\bigskip

\noindent [19] Garban, C. $\&$ Sepulveda, A. Statistical reconstruction of the Gaussian free field and KT transition. arXiv: 2022.122842v2 (2020).

\bigskip

\noindent [20] Gao, Y. $\&$ Zhang, F. On the chemical distance exponent for the two-sided level-set of the 2D Gaussian free field. arXiv: 2011.04955v1 (2020).

\bigskip

\noindent [21] Jerison, D., Levine, L., Sheffield, S. Internal DLA and the Gaussian free field. Duke Math. J. \textbf{163}, 2, 267-308 (2011).

\bigskip

\noindent [22] Kang, N-G $\&$ Makarov, N. G. Gaussian free field and conformal field theory. arXiv: 1101.1024v3 (2013).

\bigskip

\noindent [23] Kenyon, R. Dominos and the Gaussian free field. The Annals of Probability \textbf{29}(3) 1128-1137 (2001).

\bigskip

\noindent [24] Kesten, H. The Incipient Infinite Cluster in Two-Dimensional Percolation. \textit{Probability Theory and Related Fields} \textbf{73}, 369-394 (1986).

\bigskip

\noindent [25] Lupu, T. From loop clusters and random interlacements to the free field. Ann. Probab., \textbf{44}(3) 2117-2146 (2016).

\bigskip

\noindent [26] Rodriguez, P-F. A 0-1 law for the massive Gaussian free field. arXiv: 1505.08169v2 (2017).

 \bigskip

 \noindent [27] Rodriguez, P-F. Level set percolation for random interlacements and the Gaussian free field. Stochastic Processes and their Applications \textbf{124}(4) 1469-1502 (2014).

\bigskip

\noindent [28] Sheffield, S. Gaussian free field for mathematicians. arXiv: math/0312099v3 (2006).

\bigskip

\noindent [29] Sznitman, A. Disconnection and level-set percolation for the Gaussian free field. Journal Mathematical Society Japan. \textbf{67}(4) 1801-1843 (2015). 

\bigskip

\noindent[30] Wang, M. $\&$ Wu, H. Level-lines of Gaussian free field I: zero boundary GFF. Stochastic Processes and their Applications \textbf{127}(4) 1045-1124 (2017).

\bigskip

\noindent [31] Werner, W. Topics on the two-dimensional Gaussian Free Field. Lecture Notes.

\bigskip

\noindent [32] Werner, W. $\&$ Powell, E. Lecture notes on the Gaussian free field. arXiv: 2004.04720v2 (2021).

\bigskip

\noindent [33] Chayes, L. Discontiniuty of the Spin-Wave Stiffness in the Two-Dimensional XY Model. Commun. Math. Phys. \textbf{197}, 623-640 (1998).

\end{document}